\documentclass[12pt,a4paper]{amsart}

\title[Asymptotic dimension for covers with controlled growth]{Asymptotic dimension for covers with controlled growth}
\thanks{The first author was supported by EPSRC grant EP/V027360/1 ``Coarse geometry of groups and spaces''.}
\author{David Hume}
\address{School of Mathematics, University of Bristol, Bristol, BS8 1TX.}
\email{david.hume@bristol.ac.uk}
\author{John M. Mackay}
\address{School of Mathematics, University of Bristol, Bristol, BS8 1TX.}
\email{john.mackay@bristol.ac.uk}
\author{Romain Tessera}
\address{Universit\'e de Paris, Sorbonne Universit\'e, CNRS, Institut de Math\'ematiques
	de Jussieu-Paris Rive Gauche, F-75013 Paris, France.}
\email{romatessera@gmail.com}
\date{\today}

\usepackage{amsmath,amsthm,amsfonts,graphicx,amssymb}
\usepackage{hyperref}
\usepackage{xcolor}
\usepackage[shortlabels]{enumitem}
\usepackage{tikz}

\newtheorem{thmIntro}{Theorem}[]

\numberwithin{equation}{section}

\newtheorem{proposition}[equation]{Proposition}
\newtheorem{corollary}[equation]{Corollary}
\newtheorem{lemma}[equation]{Lemma}

\theoremstyle{definition}

\newtheorem{problem}[equation]{Problem}
\newtheorem{notation}[equation]{Notation}
\newtheorem{conjecture}[equation]{Conjecture}
\newtheorem{definition}[equation]{Definition}
\newtheorem{remark}[equation]{Remark}

\newtheorem*{theorem*}{Theorem}
\newtheorem*{proposition*}{Proposition}
\newtheorem*{assump*}{Standing assumption}
\newtheorem*{remark*}{Remark}
\newtheorem*{claim*}{Claim}

\newtheoremstyle{citing}
{3pt}
{3pt}
{\itshape}
{}
{\bfseries}
{}
{.5em}
{\thmnote{#3}}

\theoremstyle{citing}
\newtheorem*{varthm}{}

\DeclareMathOperator{\Gr}{Gr}
\DeclareMathOperator{\se}{se}
\DeclareMathOperator{\rk}{rank}
\DeclareMathOperator{\hypdim}{hypdim}
\DeclareMathOperator{\cdim}{cdim}
\DeclareMathOperator{\size}{size}
\DeclareMathOperator{\cork}{corank}

\DeclareMathOperator{\poly}{poly}
\DeclareMathOperator{\diam}{diam}

\DeclareMathOperator{\length}{length}

\DeclareMathOperator{\CAT}{CAT}

\newcommand{\bdry}{\partial_\infty}
\newcommand{\maps}{\longrightarrow}

\DeclareMathOperator{\asdim}{asdim}
\newcommand{\uasdim}{\overline{\asdim}}

\newcommand{\set}[1]{\left\{#1\right\}}
\newcommand{\setcon}[2]{\left\{#1\, :\, #2\right\}}

\newcommand{\abs}[1]{\left\lvert#1\right\rvert}

\newcommand{\cV}{\mathcal{V}}
\newcommand{\cU}{\mathcal{U}}

\newcommand{\dimT}{\dim_{\mathrm{top}}}

\newcommand{\cM}{\mathcal{M}}

\newcommand{\R}{\mathbb{R}}

\newcommand{\N}{\mathbb{N}}
\newcommand{\Z}{\mathbb{Z}}
\newcommand{\HH}{\mathbb{H}}

\def\XXint#1#2#3{{\setbox0=\hbox{$#1{#2#3}{\int}$}
		\vcenter{\hbox{$#2#3$}}\kern-.5\wd0}}

\numberwithin{equation}{section}

\begin{document}

	\begin{abstract} We prove various obstructions to the existence of regular maps (or coarse embeddings) between commonly studied spaces. For instance, there is no regular map (or coarse embedding) $\HH^n\to\HH^{n-1}\times Y$ for $n\geq 3$, or $(T_3)^n \to (T_3)^{n-1}\times Y$ whenever $Y$ is a bounded degree graph with subexponential growth, where $T_3$ is the $3$-regular tree. 
		We also resolve \cite[Question~5.2]{BenSchTim-12-separation-graphs}, proving that there is no regular map $\HH^2 \to T_3 \times Y$ whenever $Y$ is a bounded degree graph with at most polynomial growth, and no quasi-isometric embedding whenever $Y$ has subexponential growth.
		Finally, we show that there is no regular map $F^n\to \Z\wr F^{n-1}$ where $F$ is the free group on two generators. 
		
		To prove these results, we introduce and study generalizations of asymptotic dimension which allow unbounded covers with controlled growth. For bounded degree graphs, these invariants are monotone with respect to regular maps (hence coarse embeddings).
	\end{abstract}
	
	\maketitle
	\section{Introduction}
	
	The notion of a {\it coarse embedding} arises naturally in geometric group theory, the most immediate example being the inclusion of closed compactly generated subgroup of a compactly generated locally compact group (where the groups are equipped with their word metrics). More generally, whenever $G$ is a compactly generated locally compact group acting properly and continuously by isometries on a metric space $X$, then for every $x\in X$, the orbit map $g\mapsto g\cdot x$ is a coarse embedding.
	{\it Quasi-isometric embeddings} are a special case of coarse embeddings, which can be thought of as being ``undistorted''. In general, these are much better understood as they carry many more invariants than general coarse embeddings: for instance a quasi-isometric embedding between hyperbolic groups induces a quasi-symmetric map between their Gromov boundaries, while a coarse embedding need not \cite{Baker-Riley}.

	Another natural class of maps are {\it regular maps} (cf.\ \cite{BenSchTim-12-separation-graphs}): recall that a map $f:X\to Y$ between metric spaces is \textit{regular} if there is some $\kappa$ such that
	\begin{itemize}
		\item $d_Y(f(x),f(x'))\leq \kappa(1+d_X(x,x'))$ for all $x,x'\in X$; and
		\item the preimage of every open ball of radius $1$ in $Y$ is contained in a union of at most $\kappa$ balls of radius $1$ in $X$.
	\end{itemize}
	Every coarse embedding between \textit{bounded geometry}\footnote{A metric space has bounded geometry if for every $R$ there is some $C_R$ such that every ball of radius $R$ is contained in a union of at most $C_R$ balls of radius $1$.} metric spaces is a regular map but the converse is false, consider for example $\Z\to\Z, n\mapsto |n|$.  Although regular maps have been studied far less than coarse embeddings, they nevertheless appear naturally: indeed, the inclusion of any subgraph in a bounded degree graph is a regular map (while it may not be a coarse embedding).

	\subsection{Non-embedding results}
	
	We will be particularly interested in the following problem. 
	
	\begin{problem}\label{prob:stableEmbedding}
		Consider spaces $X$ and $Y$ of exponential growth (typically products of symmetric spaces and regular trees) such that there is no coarse embedding/regular map from $X$ to $Y$. Prove that there is no coarse embedding/regular map from $X$ to $Y\times Z$ where $Z$ is a space with subexponential growth (e.g.\ $\R^k$).
	\end{problem}
	
	Specific instances of this problem in the literature include \cite[Question~5.2]{BenSchTim-12-separation-graphs}, which asks whether there is a regular map $\HH^2 \to T_3 \times \R$, and \cite[Question~6.2]{Bensaid}: assume $X$ and $Y$ are
	direct products of non-Euclidean symmetric spaces such that there exists a coarse embedding from $X$ to $Y\times \R^k$ for some $k$, does this imply that $\rk(X)\leq \rk(Y)$?
	
	Our first result provides a strong negative answer to Benjamini--Schramm--Tim\' ar's question. 
	
	\begin{thmIntro}\label{thm:regTH2}
		There is no regular map $\HH^2$ to $T_3\times Z$ when $Z$ has polynomial growth.
	\end{thmIntro}
	
	It is quite easy to see that  $\HH^2$ does not admit a regular map to $T_3$ (for instance using asymptotic dimension). In \cite{buyalo2007elements}, Buyalo--Schroeder rule out the existence of a quasi-isometric embedding from $\HH^2$ to $T_3\times \R^k$, but so far, no invariant can obstruct coarse embeddings, even when $k=1$ (see \cite[Question~5.2]{BenSchTim-12-separation-graphs}). Interestingly, even the analytic techniques from \cite{HumeMackTess-PprofLie} which proved to be so efficient with direct products of rank $1$ symmetric spaces (of dimension $>2$) with Euclidean spaces fail here, see \cite[Remark 6.9]{HumeMackTess-PprofLie}.

	Regarding quasi-isometric embeddings, we improve Buyalo--Sch\-roe\-der's result as follows:
	
	\begin{thmIntro}\label{thm:QITH2}
		There is no quasi-isometric embedding from $\HH^2$ to $T_3\times Z$ when $Z$ has subexponential growth.
	\end{thmIntro}

	\medskip

	We now turn to other symmetric spaces. Recall a general symmetric space can be decomposed as a direct product of symmetric spaces $X=X_c\times X_e\times X_{nc}$, where $X_c$ is compact, $X_e$ is Euclidean, and $X_{nc}$ has no non-trivial compact or Euclidean factor.
	The \textit{rank} of a symmetric space $X$ is the maximal dimension of an isometrically embedded Euclidean space. It is denoted $\rk(X)$. Define the \textit{corank} of a symmetric space $X$ by $\cork(X)=\dim X -\rk X$ \footnote{The corank should not be confused with the ``subexponential corank'' introduced by Buyalo and Schroeder, which they also denote by $\cork$ and we will denote by $\cork_{\se}$. We will discuss this notion and its relation to corank in \S \ref{SecIntro:BS}.}.

	In \cite{BuSh-subexpcork}, Buyalo--Schroeder prove that for a pair $X,Y$ of symmetric spaces with trivial compact factors, the existence of a quasi-isometric embedding $X\to Y$ implies that $\cork(X)\leq  \cork(Y)$. 
	
	We believe that the techniques\footnote{in particular  Lemmas 14.2.2 and 14.2.4} from  \cite{buyalo2007elements} can be used to show that this inequality holds for regular maps as well.
	Here we use a different method to prove the following stronger statement.
	\begin{thmIntro}\label{thm:obsreg}
		Let $X$ be a symmetric space with trivial compact factor, let $Y$ be a direct product of real hyperbolic spaces, let $\ell,\ell'\in \N$, and let $Z$ be a bounded degree graph with subexponential growth. If there is a regular map \[X\times T_3^\ell\to Y\times T_3^{\ell'}\times Z,\] then
		\[
		\cork(X)+\ell \leq \cork(Y)+\ell'.
		\]
	\end{thmIntro}
	Let us discuss a few examples. 
	The non-existence of a regular map from $\HH^{n+1}$ to $\HH^n\times \R^k$ for any $n \geq 2, k \geq 0$ follows from previous work of the authors \cite{HumeMackTess-PprofLie}, and can most likely also be deduced as a consequence of the methods in the aforementioned result of \cite{BuSh-subexpcork}.
	Now we are able to replace $\R^k$ by any graph of subexponential growth:

	\begin{corollary}\label{cor:HntimesZ}
		For $n\geq 2$ there is no regular map  from  $\HH^{n+1}$ to $\HH^n \times Z$, whenever $Z$ is a bounded degree graph with subexponential growth.
	\end{corollary}

	 Following a subtle strategy based on filling functions outlined by Gromov in \cite[$\S7.E_2$, pages 141--144]{Gro-91-asymp-inv}, Bensaid proved that the rank of a direct product of non-Euclidean symmetric spaces and Euclidean buildings is monotonous under coarse embedding \cite{Bensaid}. 
	The monotonicity of the rank under regular maps remains open. In particular, the existence of a coarse embedding from $T_3\times T_3\times T_3$ to $\HH^k\times \HH^k$ for any $k\geq 2$ is prevented by Bensaid's monotonicity of the rank, but not the existence of a regular map. On the other hand, the existence of a coarse embedding from $T_3\times T_3\times T_3$ to $T_3\times T_3\times \R^k$ for some $k$ is not ruled out by any of the previous results. It is now a special case of the following:
	
	\begin{corollary}\label{cor:TntimesZ}
		For $n\geq 2$ there is no regular map  from $(T_3)^{n}$ to $(\HH^2)^{n-1} \times Z$ (in particular to $(T_3)^{n-1} \times Z$), whenever $Z$ is a bounded degree graph with subexponential growth.
	\end{corollary}
	
	Note that Theorem \ref{thm:obsreg} does not prevent the existence of a {\it regular map} from $T_3\times T_3\times T_3$ to $\HH^3\times \HH^3$, which is still  open.

	\subsection{The invariants}
	Let us first quickly review the known geometric invariants which are monotonous under coarse embeddings/regular maps. 
	The first one is volume growth, which for instance rules out the existence of a regular map from the $3$-regular tree $T_3$ to a graph of subexponential growth. A more refined invariant introduced by Gromov is asymptotic dimension, which can be seen as a large scale version of the covering dimension. This notion has played a fundamental role in various areas of mathematics, ranging from topology to non-commutative geometry. Its monotonicity under regular maps has been proved in  \cite{BenSchTim-12-separation-graphs}. For direct products of symmetric spaces of non-compact type, the asymptotic dimension is known to coincide with the topological dimension (see \cite[Corollary 3.6]{Carlsson-Goldfarb}): this for instance rules out the existence of a regular map from $\HH^n$ to $\HH^{n-1}$. Other invariants which are monotonous under coarse embeddings between spaces satisfying additional hypotheses include (co)homological dimension \cite{Shalom-04-Large-scale-amenable, Sauer-hom, Li-dynamchar}, homological filling functions \cite{Kropholler-Pengitore} and isoperimetric profile \cite{DKLMT}. None of these can be used to prove Theorems \ref{thm:regTH2}, \ref{thm:QITH2} or Corollaries \ref{cor:HntimesZ}, \ref{cor:TntimesZ} for any choice of $Z$ which has an infinite connected component.
	
	The invariant we use to prove Theorem \ref{thm:obsreg} is a generalisation of asymptotic dimension where we replace the usual condition that elements of the cover have uniformly bounded diameter with a requirement that they have uniformly controlled growth -- in this case, uniformly bounded from above by some subexponentially growing function. For Theorem \ref{thm:regTH2}, we need a slightly more sophisticated variant of this notion. 
	
	Given a family $\mathcal{F}$ of non decreasing functions $[0,\infty)\to[0,\infty)$, and a metric space $X$, we say that the \textbf{asymptotic dimension of $X$ modulo $\mathcal{F}$} is at most $n$ if, for any $r>0$ we can find a cover $\{X_i\}$ of $X$ such that every ball of radius $r$ in $X$ intersects at most $n+1$ elements of the cover, and that the growth of sets in the cover is uniformly bounded by some $f\in\mathcal{F}$. We distinguish between two versions of uniform control. The first (weak) version only considers the growth of individual elements of the cover. The second (strong) version also requires control on iterated unions of elements which are a bounded distance apart. We denote the first by $\asdim_\mathcal{F}(X)$ and the second by $\uasdim_\mathcal{F}(X)$. Clearly, we have $\asdim_\mathcal{F}(X)\leq \uasdim_\mathcal{F}(X)$, and we shall see an example below where this inequality is strict. The definitions we give in \S \ref{sec:RelativeAsdim} are actually more general, as they allow us to consider asymptotic dimension relative to any class of families of metric spaces (see Definition \ref{def:asdim_M}).

	The main motivation for introducing these notions is their monotonicity under regular maps.
	
	\begin{proposition}\label{prop:intro-monotone-reg} Let $X$ and $Y$ be metric spaces which are quasi-isometric to bounded degree graphs. If there is a regular map $X\to Y$, then for any collection of functions $\mathcal{F}$, we have
		\[
		\asdim_\mathcal{F}(X) \leq \asdim_\mathcal{F}(Y) \quad \textrm{and} \quad \uasdim_\mathcal{F}(X) \leq \uasdim_\mathcal{F}(Y).
		\]
	\end{proposition}
	
	When $\mathcal{F}$ is the set of constant functions then both $\asdim_\mathcal{F}(X)$ and $\uasdim_\mathcal{F}(X)$ are equal to the usual asymptotic dimension. The definition presented above is the ``covering'' definition of asymptotic dimension. We prove that there is an equivalent ``colouring'' definition and an equivalent ``polyhedral'' definition\footnote{The terminology here is taken from \cite{buyalo2007elements}} in terms of $\varepsilon$-Lipschitz maps to simplicial complexes (see Propositions \ref{prop:equivdefn} and \ref{prop:asdimsegeq-corank}). We also prove a natural product formula.
	
	\begin{proposition}\label{thm:productformula} Let $X$ and $Y$ be bounded degree graphs and let $\mathcal{F},\mathcal{F}'$ be families of functions. We have
		\begin{eqnarray*}
			\asdim_{\mathcal{FF}'}(X\times Y) \leq \asdim_{\mathcal{F}}(X)\times\asdim_{\mathcal{F}'}(Y); & \textrm{and} \\
			\uasdim_{\mathcal{FF}'}(X\times Y) \leq \uasdim_{\mathcal{F}}(X)\times\uasdim_{\mathcal{F}'}(Y) &
		\end{eqnarray*}
		where $\mathcal{FF}'$ is the family of all functions $g(n)=f(n)f'(n)$ where $f\in \mathcal{F}$ and $f'\in\mathcal{F}'$.
	\end{proposition}
	
	For convenience let us define the following sets of (non-decreasing) functions $f:[0,\infty)\to[0,\infty)$:
	\begin{itemize}
		\item $\poly(d)$ is all functions bounded from above by $Cn^d+C$ for some $C>0$,
		\item $\poly$ is all functions bounded from above by $Cn^d+C$ for some $C>0$ and some $d>0$,
		\item $\se$ is all subexponential functions, i.e. $f\in\se$ if and only if
		\[
		\lim_{r\to\infty} \frac{1}{r}\log(f(r)) = 0.
		\]
	\end{itemize}
	
	Theorem \ref{thm:obsreg} is a consequence of the following:

	\begin{thmIntro}\label{thm:asdimseformula}
		Let $X$ be a symmetric space without compact factor and let $\ell\in \N$, 
		\begin{equation}\label{eq:dimformulaIntro}
			\asdim_{\se}(X\times T_3^\ell)\geq \cork(X) + \ell.
		\end{equation}
		If, in addition, $X$ is a product of Euclidean and real hyperbolic spaces then we have 
		\begin{equation}\label{eq:Rank1dimformulaIntro}
			\asdim_{\se}(X\times T_3^\ell)=\asdim_{\poly}(X\times T_3^\ell)= \cork(X) + \ell.
		\end{equation}
	\end{thmIntro}
	Note that by Proposition \ref{thm:productformula}, \eqref{eq:dimformulaIntro} and \eqref{eq:Rank1dimformulaIntro} remain valid if we take a direct product of $X\times T_3^\ell$ with a bounded degree graph of subexponential growth. 
	
	The key step for the lower bound $(\ref{eq:dimformulaIntro})$ is that any cover with uniformly subexponentially growing sets intersects an exponentially distorted subset in uniformly bounded sets. This heuristic works out to give us:
	
	\begin{proposition}\label{prop:expdist} Let $X, Y$ be uniformly discrete metric spaces. If there is a regular map $f:Y\to X$, and a constant $C$ such that for all $y,y'\in Y$
		\begin{equation}\label{eq:expdist}
			d_X(f(y),f(y')) \leq C\log(1+d_Y(y,y'))+C
		\end{equation}
		then $\asdim_{\se}(X)\geq \asdim(Y)$. 
	\end{proposition}
	To obtain the lower bound we require a variant of this result which weakens \eqref{eq:expdist} to only hold for a fixed $y$ (and all $y'$) at the expense of assuming a stronger hypothesis on $Y$ (that it is quasi-isometric to a Carnot group).

	The key point in proving the upper bound in Theorem \ref{thm:asdimseformula} is the construction of a particular cover of the hyperbolic plane (see Proposition \ref{prop:tesselationH2}), which is then extended to the hyperbolic space $\HH^d$ by considering a natural quasi-isometric  embedding of $\HH^d\to(\HH^2)^{d-1}$. A careful analysis of this embedding yields the following lower bound.
	
	\begin{thmIntro}\label{prop:hypplanecover} For all $d\geq 2$,
		\[
		\asdim_{\poly(d)}(\HH^d)\leq d-1.
		\]
	\end{thmIntro}

	The definition of $\uasdim$ is somehow designed to rule out covers like the one constructed for the proof of Theorem \ref{prop:hypplanecover}.  And indeed, we have:
	
	\begin{thmIntro}\label{thm:StrongPolyDimH2} For every $d\in\N$, we have
		\[
		\uasdim_{\poly(d)}\HH^2 = 2.
		\]
	\end{thmIntro}
	Note that Theorem \ref{thm:regTH2} follows from Theorem \ref{thm:StrongPolyDimH2} (Theorem \ref{thm:QITH2} will be a consequence of its proof).

	We show that  $\asdim_\mathcal{F}$ and $ \uasdim_\mathcal{F}$ satisfy a fibering stability property, which in particular implies the following result. 
	\begin{proposition}\label{prop:asdimext} 
		Suppose that $G$ is a finitely generated group, and that we have a short exact sequence
		\[1\to P \to G \to Q\to  1,\]
		where each finitely generated subgroup of $P$, equipped with the induced metric, has subexponential growth. 
		Then  $\uasdim_{\se}(G)\leq  \asdim(Q)$.
	\end{proposition}
	This result applies to wreath products.
	\begin{corollary}\label{cor:asdimse-wreath}
		If $L$ is a finitely generated group with subexponential growth and $Q$ is a finitely generated group with finite asymptotic dimension, then the wreath product $L \wr Q = \bigoplus L \rtimes Q$ satisfies 
		\[
			\asdim_{\se}(L\wr Q) \leq \asdim(Q).
		\]
	\end{corollary}
\begin{proof}
	Any finitely generated subgroup of $\bigoplus L$ is contained in a subgroup $H \cong L^n$ of $\bigoplus L = \{f: Q \to L\}$ where all functions have support in some common finite set of size $n$.
	Since $H \cong L^n$ has subexponential growth, and the inclusion of $H\cong L^n$ (with the $\ell^1$-metric) into $L \wr Q$ is distance-non-contracting, we have that $H$ has subexponential growth with respect to the induced metric of $L \wr Q$.  Proposition~\ref{prop:asdimext} applies to give the result.
\end{proof}
	For example, combining this corollary with (\ref{eq:dimformulaIntro}) we deduce that $F^n$ does not regularly embed into $\Z\wr F^{n-1}$, where $F$ denotes the free group on two generators (observe that the asymptotic dimension of $\Z\wr F^{n-1}$ is infinite).
	
	\subsection{Comparison with Buyalo--Schroeder's invariants}\label{SecIntro:BS}
	As we have already mentioned, the spaces considered in Theorem \ref{thm:obsreg} also feature prominently in the research of Buyalo--Schroeder, who construct several invariants designed to obstruct {\it quasi-isometric embeddings} between products of symmetric spaces of non-compact type, regular trees and Euclidean spaces. Three notable invariants from their works are the hyperbolic dimension ($\hypdim$), hyperbolic rank ($\rk_h$) and subexponential corank ($\cork_{\se}$). 
	Their notion of hyperbolic dimension was our main source of inspiration for introducing our own relative notions of asymptotic dimension. Indeed, it can be seen as asymptotic dimension modulo pieces satisfying some uniform variation on a doubling condition: i.e.\ the existence of a constant $N$ such that for all $R\geq 1$, balls of radius $2R$ can be covered by at most $N$ balls of radius $R$. Such a condition is monotonous under quasi-isometry, but not under coarse embedding. Hence the hyperbolic dimension has no reason {\it a priori} to be monotonous under coarse embedding.
	Observe that 
	\[
	\hypdim(X) \geq \asdim_{\poly}(X) \geq  \asdim_{\se}(X).
	\]
	The first inequality follows from the observation that uniformly doubling spaces have uniform polynomial growth. The second inequality holds since polynomial functions are subexponential. The following inequality is more subtle.
	
	\begin{thmIntro}\label{thm:comparedim} Let $X$ be a metric space which is quasi-isometric to a bounded degree graph. We have
		\[\asdim_{\se}(X) \geq \max\{\cork_{\se}(X),\rk_h(X)\}\]
	\end{thmIntro}
	For subexponential corank, Theorem \ref{thm:comparedim} is a natural consequence of the ``polyhedral'' definition of $\asdim_{\se}$ while for hyperbolic rank, we apply Proposition~\ref{prop:expdist}/Proposition \ref{prop:expdistlb} to a family of ``level sets'' in a hyperbolic space $Y$ which are exponentially distorted, and whose asymptotic dimension (as a family) is at least the topological dimension of the boundary of $Y$.
	We may have $\asdim_{\se}(X)>\rk_h(X)$, for example $\asdim_{\se}(T_3)=1>0=\rk_h(T_3)$.  We do not have an example where $\asdim_{\se}(X)>\cork_{\se}(X)$, though perhaps one exists, since one can show that $\asdim_{\poly(1)}(\HH^2)=2$ while a `linear corank' analogous to that of Buyalo--Schroeder has value $1$ (since the preimages of a horocyclic projection $\HH^2\to\R$ grow linearly).
	
	\subsection{Questions/Conjectures}
	\begin{conjecture} Let $X$ be a symmetric space with trivial compact factor and let $\ell\in\N$. We have
		\begin{equation}\label{eq:dimformulaIntro2}
			\asdim_{\poly}(X\times T_3^\ell)=\asdim_{\se}(X\times T_3^\ell)=\cork(X)+ \ell.
		\end{equation}
	\end{conjecture}
	
	We also conjecture that Theorem \ref{prop:hypplanecover} is optimal in the following sense.
	\begin{conjecture} 
		We have  $\asdim_{o(\poly(d))}(\HH^d)= d,$ where $o(\poly(d))$ is the family of non-decreasing functions $f$ such that $\lim_{t\to \infty}f(t)/t^d=0$.
	\end{conjecture} 
	The weaker equality $\asdim_{\poly(d-1)}(\HH^d)= d$, together with Proposition \ref{thm:productformula} would be enough to rule out the existence of a regular map $\HH^d\to\HH^{d-k}\times (T_3)^k$ whenever $1\leq k\leq d$.
	This would be optimal as there exists a quasi-isometric embedding $\HH^d\to \HH^{d-k}\times\HH^{k+1}\to\HH^{d-k}\times (T_3)^{k+1}$ for all $0\leq k\leq d-2$ \cite{BradyFarb-98-filling,BuDrSc-07-hyp-grp-prod-trees}.
	
	We believe that $\uasdim_{\poly(n)}\HH^d= d$ holds for every $d,n$, and more generally:
	
	\begin{conjecture}
		For every rank one symmetric space $X$ of noncompact type and every $d$, we have $\uasdim_{\poly(d)}(X)=\asdim(X)$.
	\end{conjecture}
	
	Finally, Theorems  \ref{thm:regTH2} and Theorem \ref{thm:QITH2} should both be special cases of the following conjecture.
	\begin{conjecture}
		There is no regular map $\HH^2\to T_3 \times Y$ for any space $Y$ of subexponential growth.
	\end{conjecture}
	
	\subsection{Plan of paper}
	In Section~\ref{sec:RelativeAsdim} we define our invariants, and show their basic properties such as their behaviour under regular maps, products and extensions (Propositions~\ref{prop:intro-monotone-reg}, \ref{thm:productformula} and \ref{prop:asdimext}).
	In Section~\ref{sec:lowerbounds-asdimse} we use exponentially distorted subsets to establish the lower bound in Theorem~\ref{thm:asdimseformula}; we also prove Theorem~\ref{thm:comparedim} relating $\asdim_{\se}$ to $\cork_{\se}$ and $\rk_h$.
	In Section~\ref{sec:upper-bounds-hyp-plane} we construct decompositions of the real hyperbolic plane (and consequently also higher dimensional real hyperbolic spaces), establishing the upper bound in Theorems~\ref{thm:asdimseformula} and \ref{prop:hypplanecover}. (Theorem~\ref{thm:obsreg} follows immediately from Theorem~\ref{thm:asdimseformula}.)
	Finally, in Section~\ref{sec:unif-asdimse-hyp-plane} we show the uniform asymptotic dimension $\uasdim_{\poly(d)} (\HH^2)=2$ (Theorem~\ref{thm:StrongPolyDimH2}), and the non-embedding results Theorems~\ref{thm:regTH2} and \ref{thm:QITH2}.


\section{Asymptotic dimension relative to classes of metric families}\label{sec:RelativeAsdim}
In this section we introduce a relative version of asymptotic dimension.  Since we will have to work with metric spaces which happen to be subsets of vertices of graphs equipped with the induced metric, it will be convenient to work with the following class of metric spaces:
a metric space $X$ is ($1$-)\textbf{uniformly discrete} if $d(x,y)\geq 1$ whenever $x\neq y$.

In the sequel it will be more convenient to work with families of metric spaces, that we shall simply refer to as ``metric families'' rather than single metric spaces. A metric family $\mathcal X$ is said to have {\it bounded geometry} if for all $R>0$, there exists $N>0$ such that balls of radius $R$ in any metric space from $\mathcal X$ have cardinality at most $N$. In the sequel all of our metric families will be $1$-uniformly discrete and with bounded geometry.

\subsection{Definition and first properties}

We recall the following terminology adapted from \cite{GueTesYu-12-Finite-decomp-cplx}. Let $X$ be a metric space. A disjoint union $A=\coprod_iA_i\subset X$ of subsets of $X$ is called \textbf{$r$-disjoint} for some $r\geq 0$ if the $A_i$ are pairwise at distance at least $r$: this will be denoted by $A=\coprod_i^{r\textrm{--disjoint}} A_i$.

\begin{definition}
Let $\mathcal X$ and $\mathcal Y$ be metric families, let $r\geq 0$, and let $d\in \N$. We write  $\mathcal X\stackrel{(r,d)}{\maps}\mathcal Y$
 if for all $X\in \mathcal X$, we can write $X=X_0\cup\cdots\cup X_d$, such that 
for each $0\leq i\leq d$, $X_i= \coprod_j^{r\textrm{--disjoint}}X_{ij}$, where $\{X_{ij}\}\subset \mathcal Y$.
	We now let $\mathcal{M}$ be a class of metric families. We say that a metric family $\mathcal X$ \textbf{$d$-decomposes over $\cM$} , 
if for every $r\geq 0$, there exists $\mathcal Y\in \mathcal M$ such that $\mathcal X\stackrel{(r,d)}{\maps}\mathcal Y$.
\end{definition}

We shall need the following reinforcement of the notion of decomposition. First we need to introduce the following notation. 
\begin{notation}\label{notation:N}
	Given a metric space $X$ and $\mathcal Y$ a family of subspaces of $X$, for all $s>0$ and $m\geq 0$, $N^m_s(\mathcal Y)=\{N^m_s(Y)\mid Y\in \mathcal Y\}$ is defined inductively as follows: for all $Y\in \mathcal Y$, $N^0_s(Y)=Y$, and for each $k \geq 0$, $N^{k+1}_s(Y)$ is the union of all the $Y'\in \mathcal Y$ which intersect the closed $s$-neighbourhood of $N^{k}_s(Y)$. 
\end{notation}

\begin{definition}
	Let $\mathcal{M}$ be a class of metric families. We say that a metric family $\mathcal X$ \textbf{uniformly $d$-decomposes over $\cM$}, 
if for every $r\geq 0$, there exists $\mathcal X\stackrel{(r,d)}{\maps}\mathcal Y$, where $\mathcal Y$ is a family of metric subspaces of $\mathcal X$ such that $N^m_s(\mathcal Y)\in \mathcal M$ for all $m,s\geq 1$. 
\end{definition}

We now introduce certain classes of families of metric spaces. 

\begin{enumerate}
	\item A metric family $\mathcal X$ is \textbf{bounded} if there exists $R\geq 0$ such that all metric spaces in $\mathcal X$ have diameter at most $R$. We denote by $\mathcal B$ the class of bounded families.
\item Given an increasing function $V:\N\to \N$, a metric family $\mathcal X$ has growth function $\lesssim V$ if there exists $C\geq 1$ such that for all $r\geq 0$, balls of radius $r$ in all metric spaces in $\mathcal X$ have cardinality at most $CV(Cr)+C$. We denote by $\mathcal M_V$ the class of metric families of growth $\lesssim V$.  
\item $\mathcal M_{\se}$ is the class of metric families of subexponential growth, i.e.\
$\mathcal M_{\se}=\bigcup_V\mathcal M_V,$
where $V$ runs through all subexponential non-decreasing functions. 
\item For every $d\geq 1$, we let $\mathcal M_{\poly(d)}=\mathcal M_{r\mapsto r^d}$.
\item Finally we let $\mathcal M_{\poly}=\bigcup_{d\geq 1}\mathcal M_{\poly(d)}$.
\end{enumerate}

We observe that a metric space $X$ satisfies  $\asdim X\leq d$ if and only if $\{X\}$ (uniformly) $d$-decomposes over $\mathcal B$. By extension, we shall say that a metric family $\mathcal X$ has asymptotic dimension at most $d$ if $\mathcal X$ $d$-decomposes over $\mathcal B$. More generally we define:

\begin{definition}\label{def:asdim_M}
Let $\mathcal M$ be a collection of metric families and let $\mathcal X$ be a metric family.
\begin{itemize}
	\item The \textbf{asymptotic dimension of $\mathcal X$ modulo $\mathcal M$}, denoted by $\asdim_{\mathcal M} \mathcal X$,  is the smallest $d\in \N\cup \{\infty\}$ such that $\mathcal X$ $d$-de\-com\-pos\-es over $\mathcal M$.
	\item The \textbf{{\it uniform} asymptotic dimension of $\mathcal X$ modulo $\mathcal M$}, denoted by $\uasdim_{\mathcal M} \mathcal X$, is the smallest $d\in \N\cup \{\infty\}$ such that $\mathcal X$ {\it uniformly} $d$-decomposes over $\mathcal M$.
\end{itemize}
\end{definition}

\begin{remark}  In this paper, we will focus on the following cases:
\begin{enumerate}
\item By definition, we have $\asdim=\asdim_\mathcal B=\uasdim_\mathcal B$.
\item To be coherent with the notation given in the introduction, we shall denote
$\asdim_{\se}=\asdim_{\mathcal M_{\se}}$, $\asdim_{\poly(d)}=\asdim_{\mathcal M_{\poly(d)}}$, and $\asdim_{\poly}=\asdim_{\mathcal M_{\poly}}$ (and similarly for $\uasdim$).
\end{enumerate}
\end{remark}
We now study the monotonicity of these notions under regular maps. A map from $\mathcal X\to \mathcal Y$ is a family of maps $f=\{f_X\mid X\in \mathcal X\}$ such that the target space of $f_X$ is a space from $\mathcal Y$. Such a map is regular if there exists $C,m$ such that the maps $f_X$ are $C$-Lipschitz and have pre-images of cardinality at most $m$. For $1$-uniformly discrete metric spaces this is easily seen to be equivalent to the definition of regular map given in the introduction.

For convenience we will simply write $f$ instead of $f_X$.
A class $\mathcal M$ will be called \textbf{stable under regular maps} if for all metric families $\mathcal X$ which admit a regular map to a metric family contained in $\mathcal M$, we have $\mathcal X\in \mathcal M$.
Note that the classes $\mathcal B$ and $\mathcal M_V$ are stable under regular maps. 

A remark that we will exploit at various places is the fact that if $X$ is a metric space (satisfying our standing assumptions of uniform discreteness and bounded geometry), and $A\subset X$, and $r\geq 0$ then the map from the closed $r$-neighbourhood of $A$ (which we denote by $[A]_r$) to $A$ which to every $x$ assigns some closest point $x'$ in $A$ is a regular map: it is $(2r+1)$-Lipschitz and its pre-images have cardinality at most the maximal size of a ball of radius $r$ in $X$.

\begin{proposition}[Stability under regular maps, Proposition~\ref{prop:intro-monotone-reg}]\label{prop:regStab}
Let $\mathcal M$ be a class of metric spaces which is stable under regular maps. Let $\mathcal X$ be a metric family which regularly maps to a metric family $\mathcal Y$. Then 
\[\asdim_{\mathcal M} \mathcal X\leq \asdim_{\mathcal M} \mathcal Y, \]
and 
\[\uasdim_{\mathcal M} \mathcal X\leq \uasdim_{\mathcal M} \mathcal Y.\]
\end{proposition}
\begin{proof}
	Denote $d=\asdim_{\mathcal M} \mathcal Y$. Then, for all $r$, there exists a decomposition $\mathcal Y\stackrel{(r,d)}{\maps}\mathcal Z$, where $\mathcal Z\in \mathcal M$. Now fix a regular map $f:\mathcal X\to \mathcal Y$. We observe that since $f$ is $C$-Lipschitz, if two sets $A,B\subset Y$, for some $Y\in \mathcal Y$, are at distance at least $r$, then their pre-images by $f$ are at distance at least $r/C$. Hence $\mathcal X\stackrel{(r/C,d)}{\maps}\{f^{-1}(Z)\mid Z\in \mathcal Z\}$.
Note that $f$ induces a regular map from $\{f^{-1}(Z)\mid Z\in \mathcal Z\}$ to $\mathcal Z$. Since $\mathcal M$ is stable under regular maps, we deduce that $\{f^{-1}(Z)\mid Z\in \mathcal Z\}\in \mathcal M$. Hence we have proved that $\asdim_{\mathcal M} \mathcal X\leq d$, which is the first statement of the proposition.
The case of uniform asymptotic dimension is similar, the only additional argument which is needed is the fact that for all $s,m$, and $Z\in \mathcal Z$, 
\[f(N^m_{s}(f^{-1}(Z)))\subset N^m_{Cs}(Z),\]
which follows once again from the fact that $f$ is $C$-Lipschitz. Indeed, this inclusion implies that $f$ induces a regular map from  $\{N^m_{s}(f^{-1}(Z))\mid Z\in \mathcal Z\}$ to $\{N^m_{Cs}(Z)\mid Z\in \mathcal Z\}$, and we conclude from the fact that $\mathcal M$ is stable under regular maps.
\end{proof}

\subsection{Fibering property}
We now work to establish a fibering property for these dimensions.  A class $\mathcal M$ of metric families is \textbf{stable under inclusion} if whenever $\mathcal X\subset \mathcal Y$ and $ \mathcal Y\in \mathcal M$, then $\mathcal X\subset \mathcal M$. Since the identity is a regular map, we deduce that a class of metric families which is stable under regular maps is also stable under inclusion.

\begin{definition}\label{def:overlineM}
Let $\mathcal M$ be a class of uniformly discrete metric families. Define the class $\mathcal D_0(\mathcal M)$ to be the class of metric families which $0$-decompose over $\cM$.
\end{definition}
 In other words, $\mathcal D_0(\mathcal M)$ consists of metric families
$\mathcal X$ such that for all $r$, there exits $\mathcal Y\in \mathcal M$ such that $\mathcal X\stackrel{(r,0)}{\maps} \mathcal Y$.

 For example, for $X = \{(s,t^2) : s, t \in \N\} \subset \Z^2$, we have $\{X\}\notin \mathcal M_{\poly(1)}$ but $\{X\} \in \mathcal D_0(\mathcal M_{\poly(1)})$.
We observe that if $\mathcal M$ is stable under inclusion (respectively regular maps), then so is $\mathcal D_0(\mathcal M)$.

\begin{lemma}\label{lem:asdim-composition}
Let $\mathcal M$ be a class of uniformly discrete metric spaces which is stable under inclusion. Then \[\asdim_{\mathcal M}=\asdim_{\mathcal D_0(\mathcal M)}\quad {\text and }\quad \uasdim_{\mathcal M}=\uasdim_{\mathcal D_0(\mathcal M)}.\]
\end{lemma}

\begin{proof}
We start with $\asdim_{\mathcal M}$. Clearly since $\mathcal M\subset \mathcal D_0(\mathcal M)$, we have that $\asdim_{\mathcal D_0(\mathcal M)}\leq \asdim_{\mathcal M}$. To prove the reverse inequality we let $\mathcal X$ be a metric family such that 
$\asdim_{\mathcal D_0(\mathcal M)} \mathcal X\leq d$, meaning that for each $X\in \mathcal X$, $X=X_0\cup\cdots \cup X_d$ such that each $X_i=\coprod_j^{r\textrm{--disjoint}}X_{ij}$ and $\{X_{ij}\mid i,j,X\}\in \mathcal D_0(\mathcal M)$. This implies that for all $i,j,X$, we have $X_{ij}=\coprod_k^{r\textrm{--disjoint}}X_{ijk}$, and $\{X_{ijk}\mid i,j,k,X\}\in \mathcal M$. In other words $\asdim_{\mathcal M}\mathcal X\leq d$ and we are done.

Let us turn to $\uasdim_{\mathcal M}$, which is a bit more subtle. Once again, we let $\mathcal X$ be such that $\uasdim_{\mathcal D_0(\mathcal M)} \mathcal X\leq d$, meaning that for each $X\in \mathcal X$, $X=X_0\cup\cdots \cup X_d$ such that each $X_i=\coprod_j^{r\textrm{--disjoint}}X_{ij}$ and such that $\{N_{s}^m(X_{ij})\mid i,j,X\}\in \mathcal D_0(\mathcal M)$ for all $s,m$. Since $\mathcal M$, and therefore $\mathcal D_0(\mathcal M)$ are stable under inclusion, we can assume without loss of generality that each $X_{ij}$ is $r$-connected. A straightforward induction on $m$ shows that $N_{s}^m(X_{ij})$ is therefore $(s+r)$-connected. Now pick $R> r+s$. Since $\{N_{s}^m(X_{ij})\mid i,j,X\}\in \mathcal D_0(\mathcal M)$, we can write
	\[N_{s}^m(X_{ij})=\coprod_k^{R\textrm{--disjoint}}Z_{ijk}\] with  $\{Z_{ijk}\mid i,j,k,X\}\in \mathcal M$. But since $N_{s}^m(X_{ij})$ is $(s+r)$-connected, our choice of $R$ forces $Z_{ijk}$ to be empty except for a single $k$. 
In other words, $N_{s}^m(X_{ij})\in \mathcal M$, and so we are done. \end{proof}

\begin{proposition}[Fibering property]\label{prop:fibering}
	Let $\mathcal M$ be a class of uniformly discrete metric spaces which is stable under inclusion. Let $f:\mathcal X\to \mathcal Y$ be a Lipschitz map such that for all $R$, the family \[\{f^{-1}(A)\mid  A\subset Y\in \mathcal Y, \diam(A)\leq R \}\in \mathcal D_0(\mathcal M).\]
Then $\uasdim_{\mathcal M}(\mathcal X)\leq \asdim(\mathcal Y).$ \end{proposition}\label{prop:Lfibering}
\begin{proof}
Let $C$ be the Lipschitz constant of $f$. 
Let $d=\asdim(\mathcal Y)$, which we assume to be finite. Then, for all $r$, $\mathcal Y\stackrel{(r,d)}{\maps}\mathcal Z$, where $\mathcal Z$ is a bounded family of subspaces. As in the proof of Proposition \ref{prop:regStab}, we have that 
 $\mathcal X\stackrel{(r/C,d)}{\maps}\{f^{-1}(Z)\mid Z\in \mathcal Z\}$. Again like in the proof of that proposition, we have 
 \[N^m_{s}(f^{-1}(Z))\subset f^{-1}(N^m_{Cs}(Z)),\]
for all $Z\in \mathcal Z$. Since $\mathcal Z$ is a bounded family,  so is $\{N^m_{Cs}(Z)\mid Z\in \mathcal Z\}$: let $R$ be an upper bound on the diameter of its subsets. We have
\[\{N^m_{s}(f^{-1}(Z))\mid Z\in \mathcal Z\}\subset \{f^{-1}(A)\mid  A\subset Y\in \mathcal Y, \diam(A)\leq R \}.\]
We conclude that $\{N^m_{s}(f^{-1}(Z))\mid Z\in \mathcal Z\} \in \mathcal D_0(\mathcal M)$, since $\mathcal D_0(\mathcal M)$ is stable under inclusion (because $\mathcal M$ is): this shows that $\uasdim_{\mathcal D_0(\mathcal M)}(X)\leq d$, and we conclude by Lemma \ref{lem:asdim-composition}.
\end{proof}

\subsection{Upper bound for group extensions}

Here we prove Proposition \ref{prop:asdimext}, which is a special case (for $\mathcal M=\mathcal M_{\se}$) of the following stability result. By abuse of notation, we say that a metric space $X\in \mathcal M$ if the family  $\{X\}\in \mathcal M$.

\begin{proposition}\label{propInSection:asdimext} 
Suppose $\mathcal M$ is stable under regular maps. Let 
\[1\to P \to G \to Q\to  1\]
be a short exact sequence of countable groups.
	Assume $G$ is equipped with a uniformly discrete proper left-invariant metric, and that each finitely generated subgroup of $P$, equipped with the induced metric belongs to $\mathcal M$, and that $Q$ is equipped with the quotient metric. Then $\uasdim_{\mathcal M}(G)\leq  \asdim(Q)$.
\end{proposition}

\begin{proof}
	Note that the subgroup $P_r$ of $P$ generated by the ball of radius $r$ in $P$ is finitely generated (because the metric is proper and so the ball is finite), and such that $P$ is an $r$-disjoint union of the left-cosets of $P_r$. 
Therefore the assumption on $P$ ensures that $P\in \mathcal D_0(\mathcal M)$ with the notation of Lemma \ref{lem:asdim-composition}.
By the choice of metrics on $G$ and $Q$, the projection $f: G\to Q$ is $1$-Lipschitz. Given $R\geq 0$, the pre-image of $B_Q(1,R)$ is the $R$-neighbourhood of $P$:
\[f^{-1}(B_Q(1,R))=[P]_R.\]
	By the previous discussion, $P\in \mathcal D_0(\mathcal M)$. Therefore, since $[P]_R$ regularly maps to $P$ and $ \mathcal D_0(\mathcal M)$ is stable under regular maps (as $\mathcal M$ is), we have that $\{f^{-1}(B_Q(q,R))\mid q\in Q\}\in \mathcal D_0(\mathcal M)$. Hence the assumptions of Proposition \ref{prop:fibering} are satisfied and we conclude that $\uasdim_{\mathcal M}(G)\leq  \asdim(Q)$.
\end{proof}

\subsection{Behaviour under direct products}
We start with a useful characterization of our relative notion of asymptotic dimension, which is classical for asymptotic dimension. The argument is not new and a version can be found in \cite{BDLM:Hurewicz}, we include a full proof here for completeness.

\begin{lemma}[Kolmogorov trick]\label{lem:asdimKol} Let $\mathcal M$ be a class of metric families which is stable under regular maps, and let $\mathcal X$ be a metric family. The following are equivalent:
	\begin{itemize}
		\item[(i)] $\asdim_{\mathcal M}(\mathcal X)\leq n$,
		\item[(ii)] for every $k\geq n$ and every $r>0$, $\mathcal X\stackrel{(r,k)}{\maps}\mathcal Z$ where $\mathcal Z\in \mathcal M$ and for all $X\in \mathcal X$, and all $x\in X$, $x$ belongs to at least $(k+1)-n$ elements of $\mathcal Z$.
		\end{itemize}
Also, the following are equivalent:
\begin{itemize}
\item[(iii)] $\asdim_{\mathcal M}(\mathcal X)\leq n$,		
\item[(iv)] for every $k\geq n$ and every $r>0$, $\mathcal X\stackrel{(r,k)}{\maps}\mathcal Z$ where for all $m,s\geq 0$, $N^m_s(\mathcal Z)\in \mathcal M$ and for all $X\in \mathcal X$, and all $x\in X$, $x$ belongs to at least $(k+1)-n$ elements of $\mathcal Z$.
\end{itemize}
\end{lemma}
\begin{proof}
This lemma is proved identically to its version for asymptotic dimension. For the convenience of the reader, we sketch it, in the first case only, as the second case is similar.
It is clear that $(i)$ and $(iii)$ are precisely the $k=n$ cases of $(ii)$ and $(iv)$ respectively, so $(ii)\Rightarrow(i)$ and $(iv)\Rightarrow(iii)$. We now prove $(i)\Rightarrow(ii)$.
We reason by induction on $k$. Let $n=\asdim_{\mathcal M}\mathcal X$. For all $r$ and all $X\in \mathcal X$, we can write $X=\bigcup_{i=0}^kX_i$ such that $X_i= \coprod_j^{3r\textrm{--disjoint}}X_{ij}$, where $\{X_{ij}\}\in \mathcal M$, and for all $x\in X\in \mathcal X$, $x$ belongs to at least $k+1-n$ distinct $X_i$. We now let $X'_{ij}=[X_{ij}]_r$. We note that $X'_i= \coprod_j^{r\textrm{--disjoint}}X'_{ij}$. We now add an other subset $X'_{k+1}$ which is the union $\bigcup_{S} Y_{S}$ where $S$ runs over subsets of $\set{1,\ldots,k}$ containing exactly $(k+1)-n$ elements, and where
\[Y_{S}=\bigcap_{s\in S}\left(X_{s}\setminus \bigcup_{i\not\in S} [X_i]_r\right).\]
We first show that the $Y_{S}$ are $r$-separated.
Consider the case where $S\neq S'$. Let $x\in Y_{S}$  and $y\in Y_{S'}$. There exists $s\notin S'$ such that $x\in 
X_{s}\setminus \bigcup_{i\not\in S} [X_i]_r$. Now $y$ does not belong to $ [X_s]_r$, so $d(x,y)>r$. 
On the other hand each $Y_S$ is contained in $X_s$ for every $s\in S$, and so is a $3r$-disjoint union of subsets of the $X_{sj}$: we let $\mathcal Y$ be the metric family consisting of these pieces. 

This gives $\mathcal X\stackrel{(r,k)}{\maps}\mathcal Z$, where $\mathcal Z=\mathcal Y\cup \{X'_{ij}\mid i,j\}$. 
Note that each element of $\mathcal Z$ is contained in a single $X'_{ij}$ for some $i\leq k$. Since $\mathcal M$ is stable under regular maps, we therefore have $\mathcal Z\in\mathcal M$.

So we are left to proving that every element $x$ is contained in at least $(k+2)-n$ distinct $X'_i$. Assume that it is contained in at most $(k+1)-n$ distinct $X'_i$. Since $X_i\subset X'_i$ for all $i\leq k$, this means that $x$ belongs to exactly $(k+1)-n$ distinct $X_i$'s: let $S$ be the subset of $\set{1,\ldots,k}$ of cardinality $(k+1)-n$ containing these indices. We claim that $x$ must belong to $Y_{S}$, hence to $X'_{k+1}$. Indeed, it belongs to $\cap_{s\in S} X_{s}$, and not to  $\bigcup_{i\not\in S} X_i'$, which is nothing but $\bigcup_{i\not\in S} [X_i]_r$. Hence it does belong to $Y_S$, and we are done.
\end{proof}

Given two classes $\mathcal M_1$ and $\mathcal M_2$ of metric families,  we denote by $\mathcal M_1 \otimes \mathcal M_2$ the class of metric families $\mathcal Y$ admitting a regular map to $\{X_1 \times X_2\mid X_i\in \mathcal X_i\}$, for some $(\mathcal X_1,\mathcal X_2)\in \mathcal M_1\times \mathcal M_2$. (We use the $\ell^1$-product metric throughout.) By construction, the class $\mathcal M_1 \otimes \mathcal M_2$ is stable under regular maps. Note that given two non-decreasing functions $V_1$ and $V_2$, we have $\mathcal M_{V_1}\otimes \mathcal M_{V_2}\subset \mathcal M_{V_1 V_2}$, where $V_1V_2$ denotes the function $t\mapsto V_1(t)V_2(t)$.

\begin{proposition}[Product formula, Proposition~\ref{thm:productformula}]\label{prop:products} Let $\mathcal M$ and $\mathcal N$ be two classes of metric families which are stable under regular maps, and let $\mathcal X$ and $\mathcal Y$ be metric families.  Then
	\begin{align*}
		\asdim_{\mathcal M \otimes \mathcal N}\{X \times Y\mid (X,Y)\in \mathcal X\times \mathcal Y\}& \leq \asdim_{\mathcal M}(\mathcal X)+\asdim_{\mathcal N}(\mathcal Y)
		\\ \intertext{ and } 
		 \uasdim_{\mathcal M \otimes \mathcal N}\{X \times Y\mid (X,Y)\in \mathcal X\times \mathcal Y\}& \leq \uasdim_{\mathcal M}(\mathcal X)+\uasdim_{\mathcal N}(\mathcal Y).
	\end{align*}
\end{proposition}
\begin{proof}
	The proof, based on Lemma \ref{lem:asdimKol}, is essentially the same as the usual product formula for asymptotic dimension.  Once again, we focus on the first case, the second one being treated similarly.
Let $m=\asdim_{\mathcal M}\mathcal X$ and $n=\asdim_{\mathcal N}\mathcal Y$, and let $k=m+n$. By Lemma \ref{lem:asdimKol}, for all $r$ and all $X\in \mathcal X$, we can write $X=\bigcup_{i=0}^kX_i$ such that $X_i= \coprod_j^{r\textrm{--disjoint}}X_{ij}$, where $\{X_{ij}\}\in \mathcal M$, and for all $x\in X\in \mathcal X$, $x$ belongs to at least $k+1-m$ distinct $X_i$ (and similarly for $\mathcal Y$).
Now we define for each $i=0\ldots, k$, $Z_i=X_i\times Y_i$. Obviously, $Z_i$ is an $r$-disjoint union of the $X_{ij}\times Y_{ik}$. Hence if we can prove that $X\times Y=\bigcup_{i=0}^kZ_i$, we will be done.
Let $(x,y)\in X\times Y$, the number of $X_i$ that do not contain $x$ is at most $m$, and the number of $Y_i$ not containing $y$ is at most $n$. Since $n+m=k$, and there are $k+1$ $X_i$'s and $Y_i$'s, this implies that $(x,y)$ must belong to at least one $X_i\times Y_i$, so we are done.  
\end{proof}

\subsection{Equivalent formulation in terms of covers}

Asymptotic dimension has several equivalent formulations. Here we give one that generalizes easily, and will be required in $\S\ref{sec:expdistlb}$. 
We recall that a cover $\mathcal U$ of a metric family $\mathcal X$ has \textbf{$R$-multiplicity at most $k$} if every closed ball of radius $R$ in each $X\in\mathcal X$ intersects at most $k$ sets in $\mathcal U$. We usually refer to $0$-multiplicity as just multiplicity. A subset $U \subset X$ is \textbf{$R$-connected} if for any $x,y \in U$ there exists $x=x_0, x_1, \ldots, x_n =y$ in $U$ with $d(x_i,x_{i+1}) \leq R$ for each $i=0, \ldots, n-1$. Note that if $\mathcal U$ has $R$-multiplicity at most $k$, then the set of all $2R$-connected components of elements of $\mathcal U$ is also a cover with $R$-multiplicity at most $k$.

\begin{proposition}\label{prop:equivdefn} Let $\mathcal M$ be a class of metric families which is stable under inclusion, and let $\mathcal X$ be a metric family. The following are equivalent.
\begin{itemize}
		\item[(i)] $\asdim_{\mathcal M}(\mathcal X)\leq n$;
		\item[(ii)] for every $R$, there is a cover $\mathcal U$ of $\mathcal X$ with $R$-multiplicity at most $n+1$, where every element of $\mathcal U$ is $2R$-connected and
		$\mathcal U \in \mathcal M.$
	
	\end{itemize}
	Similarly, the following are equivalent:
		\begin{itemize}
			\item[(iii)] $\uasdim_{\mathcal M}(\mathcal X)\leq n$;
			\item[(iv)] for every $R$, there is a cover $\mathcal U$ of $\mathcal X$ with $R$-multiplicity at most $n+1$, where every element of $\mathcal U$ is $2R$-connected and for every $s>0$ and $m\geq 0$,
			$\{N^m_s(U)\mid U\in \mathcal U\} \in \mathcal M.$	
	\end{itemize}
\end{proposition}
\begin{proof} Once again, the proof is no different than the one for asymptotic dimension. Moreover, the proof of the equivalence between $(i)$ and $(ii)$ is identical to that of $(iii)$ and $(iv)$, so we shall only prove the first one. 
	
	{\noindent $(i)\Rightarrow (ii)$:} For each $r>0$, we have $\mathcal X\stackrel{(r,n)}{\maps}\mathcal Z=\set{X_{ij}}$ with $\mathcal Z\in\mathcal M$. Set $r=2R$. Suppose $X_{ij}$ and $X_{i'j'}$ intersect some closed ball of radius $R$ in $X$. It follows that $d_X(X_{ij},X_{i'j'})\leq 2R =r$, so if $i=i'$, then $j=j'$. Thus, $\set{X_{ij}}$ has $R$-multiplicity at most $n+1$. As above, we may decompose each $X_{ij}$ into its $2R$-connected components.

\smallskip{\noindent $(ii)\Rightarrow (i)$:} Let $\mathcal V$ be a cover of $\mathcal X$ with $2R$-multiplicity $\leq n+1$ where $\mathcal V\in \mathcal M$.
Let $\mathcal X_0$ be a maximal $R$-separated collection of elements of $\mathcal V$, then for each $i=1,\ldots,n$, let $\mathcal U_i$ be a maximal collection of $R$-separated elements of $\mathcal V\setminus (\mathcal U_0\cup \cdots \cup \mathcal U_{i-1})$. 
	
Suppose $x\in X\in \mathcal X$ is not covered by $\mathcal U_0\cup \cdots \cup \mathcal U_{n}$ and choose $V\in\mathcal V$ such that $x\in V$. Since $B(x,2R)$ intersects at most $n+1$ sets in $\mathcal V$, and $V$ is not in $\mathcal U_0\cup \cdots \cup \mathcal U_{n}$,  $B(x,2R)$ intersects at most $n$ sets in $\mathcal U_0\cup \cdots \cup \mathcal U_{n}$. Therefore,
	there is some $i$ such that $B(x,2R)$ has trivial intersection with every $U\in\mathcal{U}_i$. Add $V\cap B(x,R)$ to $\mathcal U_i$. Repeating this process yields a decomposition $\mathcal X\stackrel{(R,n)}{\maps}\mathcal U=\mathcal U_0\cup\cdots\cup \mathcal U_n$. By construction, its pieces are subspaces of the pieces of $\mathcal V$, and since $\mathcal M$ is stable under inclusion, we have $\mathcal U\in \mathcal M$. Hence  $\asdim_{\mathcal M}(\mathcal X)\leq n$ and the proof is complete.
\end{proof}

\section{Lower bounds on $\asdim_{\se}$}
\label{sec:lowerbounds-asdimse}

In this section we focus on the special case where $\mathcal M=\mathcal M_{\se}$ consists of metric families of subexponential growth. Note that this class is stable under regular maps, and under direct products, i.e.\ $\mathcal M_{\se}\otimes \mathcal M_{\se}=\mathcal M_{\se}$. Recall that we denote $\asdim_{\se}=\asdim_{\mathcal M_{\se}}$.  The key goals of the section are to find lower bounds for $\asdim_{\se}$ in terms of exponentially distorted subsets ($\S\ref{sec:expdistlb}$), and show it is at least both $\rk_h$, the hyperbolic rank ($\S\ref{sec:hyprank}$), and $\cork$, the subexponential corank ($\S\ref{sec:subexpcorank}$).  

\subsection{Exponentially distorted subsets}\label{sec:expdistlb}
The proof of the lower bound of Theorem \ref{thm:asdimseformula} has three parts. 
Our first result is a kind of fibering property for $\asdim_{\se}$ under regular maps with exponential distortion. 

\begin{proposition}\label{prop:expdistlb} Let $\mathcal X,\mathcal Y$ be uniformly discrete metric families with bounded geometry, and assume $\asdim_{\se}(\mathcal X)\leq d$. Suppose there is a regular map $f=\{f_Y\}_{Y\in\mathcal Y}:\mathcal Y\to \mathcal X$, a constant $C$ and, for each $Y\in\mathcal Y$, a point $a\in Y$ such that for all $b\in Y$,
	\begin{equation}\label{eq:logdist+error}
		d_X(f_Y(a),f_Y(b))  \leq C \big(\log(1+d_Y(a,b)\big)+C.
	\end{equation}
	Then for every $r>0$, there is a decomposition $\mathcal Y\stackrel{(r,d)}{\maps}\bigcup_{Y\in\mathcal Y} \mathcal U_Y$ where each $\mathcal U_Y$ is a metric family of subsets of $Y$ satisfying the following property:
	\begin{equation}\label{eq:radsublin}
		\lim_{m\to \infty}\frac1m\max\setcon{\diam(U)}{Y\in\mathcal Y,\ U\in\mathcal U_Y,\ d_Y(a,U)\leq m} = 0.
	\end{equation}
\end{proposition} 
We say a cover is \textbf{radially sublinear} if it satisfies $(\ref{eq:radsublin})$ for some choices of $a\in Y$.

For the second part we require a definition:

\begin{definition}
	A metric space $(X,d)$ \textbf{admits dilations} if, for some $\lambda\neq 1$, there is a bijection $\psi_\lambda:X\to X$ such that for all $x,x'\in X$, we have
	\[
	d_X(\psi_\lambda(x),\psi_\lambda(x'))=\lambda d_X(x,x').
	\]
\end{definition}
Euclidean spaces $\R^n$, and more generally, Carnot groups (equipped with Carnot--Carath\'eodory metrics) are examples of spaces which admit dilations. We denote the \textbf{topological (covering) dimension} of a metric space by $\dimT$, and recall that this is the minimal $d$ such that every open covering has a refinement with multiplicity at most $d + 1$. 
Recall that a metric space is proper if its closed balls are compact.

\begin{proposition}\label{prop:dilsublintobdd} Let $N$ be a proper metric space which admits dilations and let $N'$ be a uniformly discrete metric space which is quasi-isometric to $N$. If, for all $R>0$, there is a radially sublinear cover of $N'$ with $R$-multiplicity at most $k$, then $k\geq \dimT(N)+1$.
\end{proposition}

Thirdly, we show that products of symmetric spaces and trees admit suitable embeddings of Carnot groups.
\begin{proposition}\label{prop:embeddeddil} Let $X = \prod_{i=1}^k X_i \times (T_3)^{\ell}$ where the $X_i$ are symmetric spaces of non-compact type and $T_3$ is the infinite $3$-regular tree and let $X'$ be a maximal $1$-separated subset of $X$. There is a regular map $N'\to X'$ satisfying \eqref{eq:logdist+error} where $N'$ is a uniformly discrete metric space which is quasi-isometric to a Carnot group $N$ with
	\[
	\asdim(N)=\dimT(N)=\cork(X)+\ell.
	\]
\end{proposition}
The equality $\asdim(N)=\dimT(N)$ was proved by Carlsson--Goldfarb for all simply-connected nilpotent Lie groups $N$ (\cite[Theorem 3.5]{Carlsson-Goldfarb}). Carnot groups are examples of simply-connected nilpotent Lie groups.

Before proving the three propositions above, let us explain how they combine to complete the proof of the lower bound.

\begin{proof}[Proof of Theorem~\ref{thm:asdimseformula}, lower bound (\ref{eq:dimformulaIntro})]
	Let $X$ be a symmetric space with no compact factor, and let $X'$ be a maximal $1$-net in $X\times T_3^\ell$. Fix some $R\geq 0$. 
By Proposition~\ref{prop:embeddeddil}, there exists a Carnot group $N$ so that
	\[ \asdim(N)=\dimT(N)=\cork(X)+\ell,\]
	and a regular map $N'\to X'$ satisfying \eqref{eq:logdist+error} where $N'$ is a uniformly discrete metric space which is quasi-isometric to $N$. We then deduce from Proposition \ref{prop:expdistlb} that there is a decomposition $N'\stackrel{(r,\asdim_{\se}(X))}{\maps}\mathcal N'$ where $\mathcal N'$ is radially sublinear. We deduce from the proof of Proposition \ref{prop:equivdefn} $(i)\Rightarrow (ii)$, that on refining the cover (which does not affect (\ref{eq:radsublin})), we can assume that its $R/2$-multiplicity is at most $\asdim_{\se}(X)+1$. Therefore Proposition \ref{prop:dilsublintobdd} implies that $\asdim_{\se}(X)+1\geq \cork(X)+\ell+1$, and the proof is complete.
\end{proof}

\subsubsection{Proof of Proposition \ref{prop:expdistlb}}
Proposition \ref{prop:expdistlb} is the first case of the following result. The second case (Proposition~\ref{prop:expdist}) shows that under slightly stronger hypotheses on the distortion, we can deduce that $\asdim_{\se}(X)\geq\asdim(Y)$ without assumptions on $Y$.
\begin{proposition}\label{prop:expdistlb2} Let $\mathcal X, \mathcal Y$ be uniformly discrete metric families, where $\asdim_{\se}(\mathcal X)=k$ and let $f=\{f_Y\}_{Y\in\mathcal Y}:\mathcal Y\to \mathcal X$ be a regular map. 
	\begin{enumerate}
		\item[\textup{1)}] If there is a constant $C$ such that for all $Y\in \mathcal Y$ with $f_Y:Y\to X$ there exists $a \in Y$ such that for all $b\in Y$ we have
		\begin{equation}\label{eq:logplusnorm}
			d_X(f_Y(a),f_Y(b))\leq C\log(1+d_Y(a,b))+C,
		\end{equation}
		then for every $r>0$ there is a decomposition $\mathcal Y \stackrel{(r,k)}{\maps} \mathcal V$, where $\mathcal V$ is radially sublinear.
	\item[\textup{2)}] If there is a constant $C$ such that for all $Y\in \mathcal Y$ with $f_Y:Y\to X$ and  \eqref{eq:logplusnorm} holds for all $a,b\in Y$ then $\asdim(\mathcal Y)\leq k$.
	\end{enumerate} 
\end{proposition}
\begin{proof}
Let $\lambda$ be equal to the maximum of the Lipschitz constant of $\phi$ and the maximal cardinality of its pre-images.
	Set $k=\asdim_{\se}(\mathcal X)$, fix $R>0$ and let $\mathcal X\stackrel{(\lambda R,k)}{\maps} \mathcal U$ be a decomposition with $\mathcal U\in \mathcal M_{\se}$. Then $\mathcal Y\stackrel{( R,k)}{\maps}\mathcal V$, where $\mathcal V=\setcon{f^{-1}(U)}{U\in\mathcal U}$. Note that $\mathcal V\in \mathcal M_{\se}.$  On refining $\mathcal V$, we may assume in addition that its pieces are $R$-connected. We now proceed to bound the diameter of pieces of $\mathcal V$. Note that all the pieces are finite. Indeed, assume $V$ is infinite: since it is $R$-connected, both 1) and 2) would imply that $|V\cap B(y_0,n)|$ grows exponentially with $n$. The point now is to make this quantitative.
	
	Let $V\subset Y$ be a piece with $f(V) \subset U$, and suppose $a\in Y$.
	Since $V$ is $R$-connected, $V$ contains at least $\diam(V)/R$ points, so $f(V) \subset U$ contains at least $\diam(V)/\lambda R$ points.
	On the other hand, $f(V)$ is contained in the ball about $f(a)$ of radius
	\[
		\leq C \log(1+d(a,V)+\diam(V))+C,
	\]
	so the subexponential growth of $U$ gives that
	\begin{equation}\label{eq:expdist-vol}
		\frac{\diam(V)}{\lambda R} = o\left( \log(1+d(a,V)+\diam(V)) \right).
	\end{equation}

	In case 2), we may assume $a \in V$, thus \eqref{eq:expdist-vol} implies $\diam(V) = o(\log(1+\diam(V)))$, hence $\diam(V)$ is uniformly bounded.

	In case 1), where $a$ is fixed, if radial sublinearity fails, there exists $\epsilon>0$ and a sequence of subsets $V_i$ of this form with $\diam(V_i)\to \infty$ and $\diam(V_i) \geq \epsilon d(a,V_i)$. 
	But then \eqref{eq:expdist-vol} gives $\diam(V_i) = o( \log(1+(1+\epsilon^{-1})\diam(V_i)))$, contradicting $\diam(V_i)\to\infty$.
\end{proof}

 Note that by Proposition \ref{prop:equivdefn}$(i)\Rightarrow(ii)$) under the hypothesis of 2) we can deduce that for every $R>0$ there is a radially sublinear cover $\mathcal V$ of $Y$ with $R$-multiplicity at most $\asdim_{\se}(X)+1$.

\subsubsection{Proof of Proposition \ref{prop:dilsublintobdd}}

The proposition is an immediate consequence of the following two lemmas.
\begin{lemma} Let $M$ be a metric space and let $Z$ be a maximal $1$-separated subset of $M$. If, for every $R>0$, $Z$ admits a radially sublinear cover with $R$-multiplicity at most $K$, then $M$ admits a radially sublinear open cover $\mathcal U$ with ($0$-)multiplicity at most $K$.	
\end{lemma}
\begin{proof}
	Fix a radially sublinear cover $\mathcal U$ of $Z$ with $1$-multiplicity at most $K$, so there is some $z_0\in Z$ such that
	\[
	\lim_{m\to \infty}\frac1m\max\setcon{\diam(U)}{U\in\mathcal U,\ d(U,z_0)\leq m} = 0.
	\]
	Define a cover $\mathcal V=\setcon{N^\circ_1(U)}{U\in\mathcal U}$, where $N^\circ_1(\cdot)$ denotes the open $1$-neighbourhood in $M$.
	
	Since $Z$ is a maximal $1$-separated subset of $M$, it follows that $\mathcal V$ is an open cover of $M$. We claim that it is also radially sublinear and has multiplicity at most $K$.
	
	Let $U\in\mathcal U$. We have $\diam(N^\circ_1(U))\leq \diam(U)+2$.	Therefore
	\begin{multline*}
		\max\setcon{\diam(V)}{V\in\mathcal V,\ V\cap B_M(z_0,m)\neq\emptyset} \\ \leq \max\setcon{\diam(U)+2}{U\in\mathcal U,\ U\cap B_Z(z_0,m)\neq\emptyset}. 
	\end{multline*}
	Thus $\lim_{m\to\infty} \frac1m\max\setcon{\diam(V)}{V\in\mathcal V,\ d_M(V,m_0)\leq m}=0$, as required.
	
	Secondly, we prove that $\mathcal V$ has multiplicity at most $K$. Let $m\in M$ and suppose that $m\in V_1\cap\ldots\cap V_{l}$ with $V_i\in\mathcal V$ distinct and $V_i=N^\circ_1(U_i)$ for $U_i\in\mathcal U$. We have 
	\[d_M(m,U_i) \leq 1.\]
	As $\mathcal U$ has $1$-multiplicity at most $K$, we deduce that $l\leq K$.
\end{proof}

\begin{lemma}\label{lem:sublinasdim}
	Let $(N,d)$ be a metric space with the following three properties:
	\begin{itemize}
		\item $N$ is proper (closed balls are compact);
		\item $N$ admits dilations; 
		\item $N$ admits a radially sublinear open cover $\mathcal U$ with multiplicity at most $K+1$.
	\end{itemize} 
	Then the topological dimension of $N$ is at most $K$.
\end{lemma}
\begin{proof}
	Since $N$ admits dilations, there is some $\lambda\neq 1$ and a bijection $\psi:N\to N$ such that $d_N(\psi(n),\psi(n'))=\lambda d_N(n,n')$ for all $n,n'\in N$. Replacing $\psi$ by $\psi^{-1}$ if necessary, we may assume that $\lambda\in(0,1)$.
	
	Define $B_k$ to be the closed ball of radius $k$ in $N$. 
	Let $\mathcal B_k$ be an open cover of $B_k$. Since $B_k$ is compact, there is some $n\in\N$ such that for every $x\in B_k$, there is some $B_x\in\mathcal B_k$ such that $B_k\cap B(x,\lambda^n)$ is contained in $B_x$.
	
	Now, for each $m\in\N$, we define a cover $\mathcal V_m$ of $B_k$ as follows:
	\[
	\mathcal V_m = \setcon{\psi^{m}(U)\cap B_k}{U\in\mathcal U}.
	\]
	Since $\mathcal U$ is radially sublinear, $\lim_{m\to\infty} \max_{V\in \mathcal V_m} \diam V = 0$, so we choose $m=m(n)$ such that $\max_{V\in \mathcal V_m} \diam V \leq 1/n$. It follows that $\mathcal V_m$ is a refinement of $\mathcal B_k$ with multiplicity at most $K+1$. Therefore, the topological dimension of $B_k$ is at most $K$. Via the countable union theorem, we deduce that the topological dimension of $N$ is at most $K$.
\end{proof}

\subsubsection{Proof of Proposition \ref{prop:embeddeddil}}

We will construct separate embeddings for symmetric spaces and products of trees, and combine them using the following lemmas.

\begin{lemma}\label{lem:logdistprod} For $i=1,\ldots,k$, let $N_i,X_i$ be uniformly discrete metric spaces, and let $f_i:N_i\to X_i$ be coarse Lipschitz maps satisfying \eqref{eq:logplusnorm} with respect to fixed points $a_i \in N_i$ and constants $C_i$.

	Equip $N=\prod N_i$ and $X=\prod X_i$ with the $\ell^1$ product metrics, and set $a=(a_i)$. Then the coarse Lipschitz map $f=\bigoplus f_i:N\to X$ satisfies \eqref{eq:logplusnorm} with respect to $a$ and some constant $C$.
\end{lemma}
\begin{proof} As a product of coarse Lipschitz maps, $f$ is coarse Lipschitz.
	Let $b=(b_i)\in N$ with $b_i\in N_i$, and set $C=k\max\{C_i\}$. Without loss of generality, assume that $d_{X_1}(f_1(a_1),f_1(b_1)) \geq d_{X_i}(f_i(a_i),f_i(b_i))$ for all $i$. We have
	\begin{align*}
		d_{X}(f(a),f(b)) & =  \sum_{i=1}^k d_{X_i}(f_i(a_i),f_i(b_i)) 
		 \leq  k d_{X_1}(f_1(a_1),f_1(b_1)) \\
		& \leq  C \log(1+d_{N_1}(a_1,b_1)) +C
		\\ &  \leq  C \log(1+d_{N}(a,b)) +C. \qedhere
	\end{align*}
\end{proof}

\begin{lemma}\label{lem:logdistQIinv} Let $N,X$ be metric spaces and let $q_N:N'\to N$ and $q_X:X\to X'$ be $(L,D)$-quasi-isometric embeddings. Suppose $f:N\to X$ is a coarse Lipschitz map which satisfies \eqref{eq:logplusnorm} with respect to $C$ and $a\in N$.
	Then for any $a'\in N'$ 
	the map $f'=q_X\circ f\circ q_N:N'\to X'$ 
	is a coarse Lipschitz map satisfying \eqref{eq:logplusnorm} with respect to $a'$ and some constant $C'$.
\end{lemma}
\begin{proof}
	The composition of coarse Lipschitz maps is coarse Lipschitz.
	For any $b' \in N'$ we have
	\begin{align*}
		d_{X'}(f'(a'),f'(b')) 
		& \leq 
		L d_X\big( f(q_N(a')), f(q_N(b')) \big) + D
		\\ & 
		\leq LC  \log\big(1+d_{N}(q_N(a'),q_N(b'))\big) +LC+D
		\\ &
		\leq LC \log\big(1+Ld_{N'}(a',b')+D \big) +LC+D
		\\ & \leq 
		LC\log\big(1+d_{N'}(a',b')\big)+C',  
	\end{align*}
	for some suitable $C'$, where the last step follows from $\log(1+cx+d)\leq \log((1+x)(1+c+d)) \leq \log(1+x)+c+d$.
\end{proof}

\begin{proposition}\label{prop:iwasawa}
	Let $G$ be a connected semisimple real Lie group with Iwasawa decomposition $G=KAN$ and associated symmetric space $X=G/K$. The following hold:
	\begin{enumerate}
		\item $N$ is a Carnot group;
		\item there is a constant $C$ such that the natural map $\iota:N\to X$ satisfies $d_X(\iota(a),\iota(b))\leq C\log(1+d_N(a,b))+C$ for all $a,b\in N$;
		\item if $X$ has trivial compact factor, then $\cork(X)=\dimT(N)$.
	\end{enumerate}
\end{proposition}
\begin{proof}
	The first assertion is classical (see \cite[p 172(e)]{Folland} or \cite[Theorem VI.3.4 and Theorem VI.5.1]{Helgason}). 
	For the second assertion, by Helgason \cite[Theorem VI.3.4]{Helgason} the Iwasawa decomposition $KAN$ satisfies that $N$ is `instable' in the sense of Guivarc'h~\cite[Definition A.8]{Gui-80-LoiGrandsNombres} in the solvable group $AN$, and so by Proposition A.5 and its corollary in \cite{Gui-80-LoiGrandsNombres} the embedding of $N$ into $AN$ satisfies the required distortion bound; $AN$ and $X$ are both quasi-isometric to $G$.
	 The last assertion is an obvious consequence of the definition of corank.
\end{proof}

\begin{lemma}\label{lem:ZinT3} Let $T_3$ be the infinite $3$-regular tree. There is a $3$-regular map $W:\Z\to T_3$ so that \eqref{eq:logplusnorm} holds for $W$ with $a=0$ and some $C$.
\end{lemma}
\begin{proof}
The idea is that we compose depth-first searches of depth 1, 2, 3, and so on.  

	Let $B_k$ denote the rooted binary tree of depth $k$, and label the leaves (degree $1$ vertices) of $B_k$ by elements of $\{0,1\}^k$ and order these vertices lexicographically as follows: $0^k=l_0<0^{k-1}1=l_1 < \ldots < 1^k=l_{2^k-1}$. We define a closed walk $W_k$ on $B_k$ starting at the base vertex which is a concatenation of shortest paths $o \to l_0 \to l_1 \to \ldots \to l_{2^k-1} \to o$. Covering each edge exactly twice, this walk has length $4(2^k-1)$ and meets each vertex of $B_k$ at most $3$ times. 
	
	Now take a graph with vertex set $\Z$ and edges $vw$ when $|v-w|=1$. Now connect each vertex $k$ by an edge to the root of a copy of $B_{|k|}$, and label this root vertex $k'$. The resulting graph is a subgraph of $T_3$. Now build a bi-infinite walk $W:\Z\to T_3$ on this graph as a concatenation of the following for each $k\in\Z$: starting from $k'$ follow the walk $W_{|k|}$ around the copy of $B_{|k|}$ whose root is $k'$, then take the shortest path to $(k+1)'$. This bi-infinite walk also meets every vertex at most 3 times. 

	If $W(b)$ lies in the part of the walk attached to $k\in \Z$, since the walk attached to $k$ has depth $1+|k|$ we have 
	$d_{T_3}(W(a),W(b)) \leq |k|+1+|k|=2|k|+1$.
	The case $k=0$ for \eqref{eq:logplusnorm} is trivial with $a=0,C=2$.
	Assuming $k\neq 0$, the walk from $W(a)$ to $W(b)$ has passed through a walk of depth $|k|-1$ attached to the vertex of $\Z$ one closer to $0$ than $k$, which has length $4(2^{|k|-1}-1)+2$.
	Thus \eqref{eq:logplusnorm} holds for suitable $C$:
	\begin{align*}
		d_{T_3}(W(a),W(b)) 
		& \leq 2|k|+1
		\\ & \leq C \log\big(1+ 4(2^{|k|-1}-1)\big)
		\leq C \log(1+|b|). \qedhere
	\end{align*}
\end{proof}

\begin{proof}[Proof of Proposition \ref{prop:embeddeddil}] Let $X=\prod_{i=1}^k X_i \times (T_3)^l$ where each $X_i$ is a symmetric space of non-compact type. Let $G_i$ be the connected semisimple Lie group associated to the symmetric space $X_i$, let $K_iA_iN_i$ be the Iwasawa decomposition of $G_i$ and let $\iota_i:N_i\to X_i$ be the inclusion described in Proposition \ref{prop:iwasawa}. Fix a maximal $1$-separated subset $X'$ of $X$ and a quasi-isometry $q:X\to X'$.
	
	Let $N$ be the Carnot group $\prod_{i=1}^d N_i\times \R^l$. Fix quasi-isometries $q_i:N_i'\to N_i$ where $N'_i$ is a maximal $1$-net in $N_i$. Let $N'$ be the $\ell^1$ product $\prod_{i=1}^d N'_i\times \Z^l$; this is a uniformly discrete metric space which is quasi-isometric to $N$.
	
	Define $f=\bigoplus_{i=1}^k q_i\oplus \textup{id}$ and $g=\bigoplus_{i=1}^k \iota_i\oplus \bigoplus_{j=1}^lW$ where $W:\Z\to T_3$ is the map constructed in Lemma \ref{lem:ZinT3}.
	
	By Lemmas \ref{lem:logdistprod} and \ref{lem:logdistQIinv} the map $r$ defined as the composition
	\[
	N'=\prod_{i=1}^k N'_i \times \Z^l \to_{f} \prod_{i=1}^k N_i \times \Z^l \to_{g} X \to_q X'
	\]
	is a coarse Lipschitz map satisfying \eqref{eq:logplusnorm}. We now prove that it is regular. Let $x'\in X'$. We have $\diam(q^{-1}(x'))\leq D$. As $W$ is $3$-regular, 
	\[
	\left|\bigoplus_{j=1}^l(W^{-1}(q^{-1}(x')\cap V\Z))\right|\leq 3^l(3^D+1)
	\]
	also, each $\iota_i$ is a coarse embedding, so there is some $D'$ (independent of $x'$) such that
	\[
	\diam\left(\bigoplus_{i=1}^k \iota_i^{-1}(q^{-1}(x')\cap \prod_{i=1}^k N_i\right)\leq D'.
	\]
	Hence $\diam(g^{-1}(q^{-1}(x')))$ is contained in a uniformly bounded number of balls of uniformly bounded diameter. As $N'$ has bounded geometry (so balls in $N'$ of any fixed radius contain a uniformly bounded number of points) and $f$ is a quasi-isometry (being a direct sum of finitely many quasi-isometries), we have a uniform upper bound on $|r^{-1}(x')|$, as required.
\end{proof}

\subsection{Hyperbolic rank}\label{sec:hyprank}

Consider the following ways to measure the hyperbolic ``size'' of a space.

\begin{definition}
	The (regular/coarse/quasi-isometric)-\textbf{hyperbolic rank} of a metric space $X$ ($\rk_h^r(X)$/$\rk_h^c(X)$/$\rk_h^q(X)=\rk_h(X)$) is the maximal $k$ such that there is a regular/coarse/quasi-isometric embedding $Y\to X$ where $Y$ is a bounded degree visual hyperbolic graph whose boundary has topological dimension $k$.
\end{definition}
For quasi-isometric embeddings, this notion was introduced by Gromov (as a `corank' \cite[\S 6.B$_2'$]{Gro-91-asymp-inv}), and developed by Buyalo--Schroeder (with the variation that they required $Y$ to be a $\CAT(-1)$ Hadamard space \cite{BuSh-subexpcork}).

We recall that a hyperbolic metric space $X$ is \textbf{visual} if there is a point $x_0\in X$ such that every $y\in X$ is within a uniform distance of a geodesic ray in $X$ starting at $x_0$. It is natural to require it of $Y$ in the definition of hyperbolic rank since we are interested in the boundary of $Y$; the two notions (with/without visual) are the same since one can replace any non-visual bounded degree hyperbolic graph by the (convex hull of) all geodesic rays from some basepoint, without changing the boundary.

The following definition is useful to study the topological dimension of the boundary of a hyperbolic group.
\begin{definition}[{Buyalo \cite[Proposition 3.2]{Bu-05-asdim-capdim}}]
	The \textbf{capacity dimension} of a metric space $Z$ is the infimum $\cdim(Z)$ of all integers $m$ with the following property: there exists a constant
	$c > 0$ such that for all sufficiently small $s > 0$, $Z$ has a $cs$-bounded covering
	with $s$-multiplicity at most $m + 1$.
\end{definition}
(This notion is also called ``linearly controlled metric dimension'', see \cite[Chapter 11]{buyalo2007elements}.)
It is always the case that $\cdim(Z)\geq\dimT(Z)$.

If $G$ is a hyperbolic group then $\asdim(G)=\cdim(\partial G)+1$ \cite{Bu-05-asdim-capdim, BuLe-selfsim}.  We now show $\asdim_{\se}$ is bounded below by the capacity dimension of the boundary.
\begin{proposition}\label{prop:hypgraphlb}
	Let $X$ be a visual hyperbolic  bounded degree graph. There is a uniformly discrete metric space $N$ with asymptotic dimension at least $\cdim(\partial_\infty X)$ and a map $r:N\to X$ satisfying \eqref{eq:logplusnorm} for all $a,b \in N$. In particular, $\asdim_{\se}(X)\geq \cdim(\partial_\infty X)$.
\end{proposition}
This has the following immediate consequence.
\begin{corollary}\label{asdimsegeqhypdim} Let $X$ be a bounded degree graph. We have \[ \asdim_{\se}(X)\geq \rk_h^r(X)\geq \rk_h^c(X)\geq \rk_h^q(X). \]
\end{corollary}
\begin{proof}
	The final two inequalities are obvious as quasi-isometries are coarse embeddings, and coarse embeddings between bounded degree graphs are regular maps. Suppose there is a regular map $Y\to X$ where $Y$ is a hyperbolic graph with bounded degree. By Proposition \ref{prop:hypgraphlb}, $\asdim_{\se}(G)\geq \cdim(\partial_\infty G)$. Since $\dimT(Z)\leq\cdim(Z)$ for any metric space $Z$, the result follows.
\end{proof}

\begin{proof}[Proof of Proposition~\ref{prop:hypgraphlb}]
	\textbf{Step 1:} Construction of $N$.
	Let $Y'$ be the hyperbolic cone of $\bdry X$ in the sense of Buyalo--Schroeder~\cite[Theorem 6.4.1]{buyalo2007elements} (cf.\ \cite{BP-03-lp-besov}).  Since $X$ is visual, $Y'$ is quasi-isometric to $X$ (both being quasi-isometric to the hyperbolic cone in the sense of Bonk--Schramm~\cite[Theorem 8.2]{BS-00-gro-hyp-embed}).

	Let us briefly recall the construction of $Y$. Rescale $\partial_\infty X$ so that it has diameter $1/2$, and denote this metric by $\rho$. For each $k\in\N$, let $X_k$ be a maximal $6^{-k}$-separated set in $\partial_\infty X$. The vertex set of $Y$ is $\bigsqcup_{k\geq 0} X_k$, with $X_0 = \{o\}$. To each $v \in X_k \subset Y$ we associate the ball $B_v:=B(v, 2\cdot 6^{-k})$. 
	Buyalo--Schroeder connect pairs of vertices $x, x' \in X_k$ if $\bar B_x \cap \bar B_{x'} \neq \emptyset$, or if $x\in X_k$ and $x' \in X_{k+1}$ and $B_{x'}\subset B_{x}$, to make a graph $Y'$.
	We add potential additional edges between $x,x' \in X_k$ if $\rho(x,x') \leq 4 \cdot 6^{-k}$ to make a graph $Y$; since such $x,x'$ have bounded distance in $Y'$, $Y$ is quasi-isometric to $Y'$ and hence $X$ too, so it suffices to find the map $r:N\to Y$.

	Define $\mathcal N=\set{(N_k,d_k)}_{k\in\N}$ to be the metric family consisting of the vertex sets of the full subgraphs $X_k$ equipped with the shortest path metric, and let $r:\mathcal N\to Y$ be the standard inclusion. 

	\textbf{Step 2:} The map $r$ satisfies \eqref{eq:logplusnorm} for all pairs of points.
	Given $x,x' \in N_k$, let $\gamma:[0,k]\to Y, \gamma':[0,k]\to Y$ be geodesics from $o$ to $x, x'$ respectively.
	By the exponential divergence of geodesics~\cite[III.H.1.26]{BH-99-Metric-spaces} there exists an exponentially growing function $e:\N\to\R$ so that if for some $t$, $d_Y(\gamma(k-t),\gamma'(k-t))\geq e(0)$, then any path joining $x$ to $x'$ in $N_k$ must have length $\geq e(t)$.
	By thinness of geodesic triangles, there exists $C$ so that if $d_Y(x,x')\geq C$ we can find $t \geq \frac{1}{2}d_Y(x,x')-C$ satisfying this condition, thus for suitable $C'$,
	\[
		d_Y(x,x') \leq C'\log(1+e(t))+C' \leq C'\log(1+d_k(x,x'))+C'.
	\]
	 
	\textbf{Step 3:} $\asdim(\mathcal N)\geq\cdim(\partial G)$. 
	For each $k$ select a partition $\setcon{U^k_{x}}{x\in X_k}$  of $\partial_{\infty}X$ so that $B(x,\tfrac{1}{2}6^{-k})\subseteq U^k_x\subseteq B(x,6^{-k})$.

	 Fix a cover $\mathcal V$ of $\mathcal N$ with $1$-multiplicity at most $\asdim(\mathcal N)+1$ where each $V\in\mathcal V$ has diameter at most $M$. Now construct the cover $\mathcal U^k$ of $\partial_{\infty}X$ consisting of the sets 
	\[
		\setcon{\bigcup_{x\in N_k\cap V} U^k_x}{V\in\mathcal V}.
	\]
	For any $U\in\mathcal U^k$, we have 
	\[
		\diam_\rho(U)\leq (\diam_{d_k}(N_k\cap V)+1)\cdot 4 \cdot 6^{-k} \leq 4\cdot 6^{-k}(M+1).
	\]
	Now let $y\in\partial_\infty X$ with $y\in U^k_{x}$ and suppose $B(y,6^{-k})$ intersects $U^k_{x_0},\ldots, U^k_{x_l}$ at $y_0,\ldots, y_l$ respectively. It follows that for each $i$
	\[
		\rho(x,x_i)\leq \rho(x,y)+\rho(y,y_i)+\rho(y_i,x_i) \leq 3 \cdot 6^{-k}.
	\]
	Hence $d_k(x,x_i)\leq 1$. Since $\mathcal V$ has $1$-multiplicity at most $\asdim(\mathcal N)+1$, $\{x_i\}$ intersects at most $\asdim(\mathcal N)+1$ of $\mathcal V$. Hence, $\mathcal U^k$ has $6^{-k}$-multiplicity at most $\asdim(\mathcal N)+1$. Thus, $\cdim(\partial_\infty X) \leq \asdim(\mathcal N)$.
\end{proof}

\subsection{Subexponential corank}\label{sec:subexpcorank}
Buyalo--Schroeder's subexponential co\-rank of a metric space $X$ is defined, roughly, as the minimal $n$ such that all metric spaces quasi-isometric to $X$ have a \emph{continuous} map to a space of topological dimension $n$ with subexponential fibres.  
Thus in this section we will have to move away from bounded degree graphs and uniformly discrete metric spaces.
First we build suitable Lipschitz maps.

\begin{proposition}\label{prop:asdimsegeq-corank} Let $X$ be a metric space with bounded geometry and has $\asdim_{\se}(X)\leq n$.
	Then there is a 
	Lipschitz map $p:X\to T$ where $T$ is a simplicial complex of dimension at most $n$ with the following property: 
	
	For each $x_0\in X$, each maximal separated net $X_\delta \subseteq X$ with a sufficiently large separation
	constant $\delta$, each $\sigma\geq\delta$ and every $\epsilon > 0$ there exists
	$S_0 = S_0(X_\delta,\sigma,\epsilon) \geq 1$ such that for every $S \geq S_0$ and every $t \in T$ we have
	\begin{equation}\label{subexpsize}
		\frac1S \log\abs{X_\delta\cap [p^{-1}(t)]_{\sigma}\cap B(x_0,S)}
		\leq \epsilon.
	\end{equation}
\end{proposition}
Before proving this proposition, we derive the bound on $\cork$.
Together with Corollary~\ref{asdimsegeqhypdim}, the following corollary proves Theorem~\ref{thm:comparedim}.
\begin{corollary}\label{cor:asdimsegeqcorankse}
	Let $X$ be a metric space which is quasi-isometric to a simplicial graph of bounded degree. Then $\asdim_{\se}(X)\geq\cork_{\se}(X)$.
\end{corollary}
\begin{proof}
	We recall that $\cork_{\se}(X)\leq n$ means that for any metric space $Y$ quasi-isometric to $X$ there is a continuous map
$g:Y \to T$, with the following properties
\begin{itemize}
	\item every compact subset of $T$ has topological dimension at most $n$,
	\item for any maximal $\delta$-separated net $Y_\delta$ of $Y$ (for sufficiently large $\delta$), for all $\sigma\geq \delta$, and all $\varepsilon>0$, there is a constant $S_0=S_0(Y_\delta,\sigma,\varepsilon)$ such that for all $S\geq S_0$ and all $t\in T$, 
	\begin{equation}\label{subexpgrowth2}
		\frac1S \log\abs{Y_\delta\cap [p^{-1}(t)]_{\sigma}\cap B(y_0,S)}
		\leq \epsilon.
	\end{equation}
	for some fixed point $y_0\in Y$.
\end{itemize}
Now suppose $\asdim_{\se}(X)\leq n$, let $Y$ be any metric space quasi-isometric to $X$ and apply Proposition \ref{prop:asdimsegeq-corank} to $Y$. We obtain a Lipschitz map $g:Y\to T$ where $T$ is a simplicial complex of dimension at most $n$ and ($\ref{subexpgrowth}$) holds. Since $g$ is Lipschitz, it is continuous. Moreover, every compact subset of $T$ has topological dimension at most $n$ and ($\ref{subexpgrowth2}$) is exactly ($\ref{subexpgrowth}$), so we deduce that $\cork_{\se}(X)\leq n$.
\end{proof}

We now return to justify the proposition.
\begin{proof}[Proof of Proposition~\ref{prop:asdimsegeq-corank}] Fix a maximal $1$-separated subset $Z\subseteq X$. 
	
	By Proposition \ref{prop:equivdefn}, for every $R>0$ there is a cover $\mathcal V\in \mathcal M_{\se}$ of $Z$ with $(R+1)$-multiplicity at most $n+1$.
	In particular, this means that for every $\varepsilon>0$ there is some $S_0$ such that for all $V\in\mathcal V$, any $v_0\in V$ and any $S\geq S_0$ we have
	\begin{equation}\label{subexpgrowth}
		\frac1S \log |V\cap B_Z(v_0,S)| <\varepsilon.
	\end{equation}
	
	Now define $\mathcal U = \setcon{N^\circ_{R+1}(V)}{V\in\mathcal V}$, where we recall that $N^\circ_{R+1}(V)$ is the open $(R+1)$-neighbourhood of $V$.

	\begin{claim*} It holds that $\mathcal U$ is a cover of $X$ with Lebesgue number at least $R$ and multiplicity at most $n+1$.
	\end{claim*}
	\begin{proof}[Proof of Claim] Recall that the Lebesgue number of a cover $\mathcal U$ is
		\[
		L(\mathcal U)=\inf\setcon{\max\setcon{d(x,X\setminus U)}{U\in\mathcal U}}{x\in X}.
		\]
		Since any $x\in X$ lies in $N^\circ_1(V)$ for some $V\in\mathcal{V}$, \[ d(x,X\setminus N^\circ_{R+1}(V))\geq R.\]
		Therefore, $L(\mathcal U)\geq R$.
		
		Suppose $x_0$ lies in $U_1,\ldots,U_k\in \mathcal U$ where $U_i=N^\circ_{R+1}(V_i)$ for some distinct $V_i\in \mathcal V$. It follows that
		$d_X(x_0,V_i)\leq R+1$. Since $\mathcal V$ has $(R+1)$-multiplicity at most $n+1$ we deduce that $k\leq n+1$. Thus $\mathcal U$ has multiplicity at most $n+1$.
	\end{proof}
	
	For each $U\in\mathcal U$, define
	\[
	\phi_U(x)=\frac{d(x,X\setminus U)}{\sum_{V\in\mathcal U} d(x,X\setminus V)}
	\]
	Note that $\phi_U$ is supported on $U$, $\sum_{V\in\mathcal U}\phi_V(x)=1$ and the image of $\phi_U$ is in $[0,1]$. Now define
	\[
	p:X\to\ell^2(\mathcal U), \quad \textup{given by} \quad p(x)(U)=\phi_U(x).
	\]
	The image of $p$ is contained in a natural realisation of the simplicial complex with vertex set $\mathcal U$ and $U_1,\ldots,U_k$ span a simplex if $\bigcap_{i=1}^k U_i\neq\emptyset$. The dimension of this simplicial complex is one less than the multiplicity of $\mathcal U$.
	
	A standard calculation shows that $p$ is Lipschitz~\cite[Lemma 9.2.2]{buyalo2007elements}.

	Finally, the pre-image of a point $p^{-1}(t)$ is contained in a union of elements of $\mathcal U$ which contain a common point. We must therefore verify that these unions satisfy $(\ref{subexpsize})$. It suffices to prove this holds for each $U\in\mathcal U$ since $\mathcal U$ has bounded multiplicity.
	
	Let $X_\delta\subset X$ be a maximal $\delta$-separated net and fix $\sigma\geq\delta$. Assume $\delta>2$, so that any map $q:X_\delta\to Z$ which maps each point in $X_\delta$ to some closest point in $Z$ is injective. Set $U=N^\circ_{R+1}(V)$ with $V\in\mathcal V$. Let $x\in X_\delta$. If $x\in [U]_\sigma$, then $q(x)$ is contained in $[U]_{\sigma+1}\subseteq [V]_{\sigma+R+2}$. Also, if $x\in B(x_0,S+\sigma)$, then $q(x)\in B(q(x_0),S+\sigma+2)$. As $q$ is injective
	\begin{align*}
		\size_{X_\delta,\sigma}(U\cap B(x_0,S)) 
		& := \abs{\left\{ x\in X_\delta : x\in [U \cap B(x_0,S)]_\sigma\right\}}
		\\
		& \leq  
		|Z\cap [V]_{sigma+R+2}\cap B(q(x_0),S+\sigma+2)| \\
		& \leq  
		M|B_{Z}(q(x_0),S+\sigma+2)\cap V|
	\end{align*}
	where $M$ is the maximal cardinality of a ball of radius $(\sigma+R+2)$ in $Z$. Fix $\varepsilon>0$. Choosing $S_0$ large enough, for all $S\geq S_0$ we have
	\begin{align*}
		& \frac1S \log\size_{X_\delta,\sigma}(U\cap B(x_0,S))  \\ 
		& \leq  
		\frac1S\log M|B_{Z}(q(x_0),S+\sigma+2)\cap V| \\ 
		& \leq  
		\frac{1}{S}\left(\log |B_{Z}(q(x_0),S+\sigma+2)\cap V| + \log M\right)  \\
		& <  \varepsilon.\qedhere
	\end{align*}
\end{proof}

\section{Upper bounds on asymptotic dimension with subexponential fibres}\label{sec:upper-bounds-hyp-plane}

The goal of this section is to prove upper bounds on $\asdim_{\se}$ for products of real hyperbolic spaces and trees (completing the proof of Theorem \ref{thm:asdimseformula}) and for group extensions. As explained in the introduction, this is achieved by constructing an explicit family of covers of $\HH^2$, and then extending to $\HH^d$, for $d\geq 2$ exploiting a certain quasi-isometric embedding $\HH^d\to (\HH^2)^{d-1}$, eventually leading to Theorem \ref{prop:hypplanecover}. 

\subsection{Upper bounds for Theorem \ref{thm:asdimseformula}}
Our first goal is to construct a special cover of $\HH^2$, leading to the following special case of  Theorem \ref{prop:hypplanecover}.

\begin{proposition}\label{prop:tesselationH2}
The hyperbolic plane $\HH^2$ satisfies $\asdim_{\poly(2)}(\HH^2) = 1$.
\end{proposition}

\begin{proof}
	The lower bound is trivial; it remains to show the upper bound. 
	Let us work in the upper half plane model. 
	In brief, the upper bound is seen by considering the tiling of $\HH^2$ partly shown in Figure~\ref{fig:tileH2}: this consists of two colours of tiles; the red $A$ tiles grow linearly (as does the blue tile containing $B_1$), while the remaining blue $B$ tiles, being quasi-isometric to $\R$ with a ray attached at each integer point, grow quadratically.

	We now provide details. 
	Let $X$ be a maximal $1$-separated subset of $\HH^2$.
	Fix $r>0$. Define $\lambda_k=\sinh kr$, for $k=0,\ldots,4$. The Hausdorff distance between the lines $x=\lambda_ky$ and $x=\lambda_ly$ is $\abs{k-l}r$. Let $S_0$ be the semicircle with boundary on the $x$-axis which has the line $x=y\sinh(4r)$ as a tangent and $(1,0)$ and $(x,0)$ as the boundary points of this semicircle with $x>1$ chosen so that $S_0$ does not cross $x=\lambda_4 y$. The dilation of factor $x^n$ for $n\in\Z$ is a hyperbolic isometry mapping $S_0$ to a new semicircle $S_n$ with boundary points $(x^{n},0)$ and $(x^{n+1},0)$. Let $S=\bigcup_{n\in\Z}S_n$.

	Let $A$ be the subset of $\HH^2$ whose boundary is the lines $x=\lambda_1y$ and $x=\lambda_3y$, let $B_1$ be the subset of $\HH^2$ whose boundary is the lines $x=0$ and $x=\lambda_1y$ and let $B_2$ be the subset of $\HH^2$ whose boundary is $S$ and the line $x=\lambda_3y$. Define $B$ to be the union of $B_2$ plus every translate of $B_1$ by an orientation-preserving isometry mapping the $y$-axis to a semicircle in $S$.
	Note that $A$ contains the $r$-neighbourhood of the line $x=\lambda_2y$ and $B$ contains the $r$-neighbourhood of $S$.
	See Figure~\ref{fig:tileH2}.
	
	\begin{figure}[h]
	\begin{tikzpicture}[yscale=0.9,xscale=0.9, vertex/.style={draw,fill,circle,inner sep=0.3mm}]
		
		\draw[black] 	(0,7)--(0,0)
		(-0.2,0)--(12,0);
		\begin{scope}
			\clip (0,0) rectangle (12,7);				
			\foreach \r in {0.015625,0.03125,0.0625,0.125,0.25,0.5,1,2,4}{
				\draw[white, very thick]
				(3*\r, 0) circle (\r cm);
				\draw[white,dotted]
				(0,0) circle (2*\r cm);
			}
		\end{scope}
		\draw[black, very thick] (0,0)--(0,7);
		
		\fill[red!70!white] (0,0) -- (4.2,7) -- (9,7) -- (0,0);
		\fill[blue!50!white] (0,0) -- (4.2,7) -- (0,7) -- (0,0);
		\begin{scope}
			\clip (0,0) rectangle (12,7);				
			\fill[blue!55!white]  (16,0) arc (0:180:4cm) arc (0:180:2cm) arc (0:180:1cm) arc (0:180:0.5cm) arc (0:180:0.25cm) arc (0:180:0.125cm) arc (0:180:0.0625cm) arc (0:180:0.03125cm) arc (0:180:0.015625cm) -- (0,0) -- (9,7) -- (16,7) -- (16,0);
		\end{scope} 
		
		\begin{scope}
			\clip (0,0) rectangle (12,7);	
			\foreach \r in {0.015625,0.03125,0.0625,0.125,0.25,0.5,1,2,4}{
				\fill[blue!50!white] (4*\r, 0) arc
				(0:180:\r cm and 0.7*\r cm) arc (180:0:\r cm and \r cm); 
				\fill[red!70!white] (4*\r, 0) arc
				(0:180:\r cm and 0.7*\r cm) arc (180:0:\r cm and 0.4*\r cm);
				\draw[black, very thick]
				(3*\r, 0) circle (\r cm);
			}
		\end{scope}
		\draw[black,very thick] (0,0)--(0,7);
		
		\begin{scope}
			\clip (0,0) rectangle (12,7);				
			\foreach \r in {0.015625,0.03125,0.0625,0.125,0.25,0.5,1,2,4}{
				\draw[black,dotted,very thick]
				(0,0) circle (2*\r cm);
			}
		\end{scope}
		\node at (4,4.7) {$A$};
		\node at (5.5,2.7) {$B_2$};
		\node at (1.7,5.8) {$B_1$};
		\node[right] at (12,4) {$C$};
		\node[above] at (4.5,7) {$x=\lambda_1y$};
		\node[above] at (9,7) {$x=\lambda_3y$};
		\node[] at (11.3,3.3) {$\psi(B_1)$};
		\node[] at (11.3,2.2) {$\psi(A)$};
	\end{tikzpicture}
		\caption{}\label{fig:tileH2}
	\end{figure}
	
	We now tesselate the right half-space of $\HH^2$ by translates of $A$ and $B$ as follows: take a semicircle $S'$ contained in the boundary of a translate of $B_2$ and an isometry $\psi$ of $\HH^2$ which maps the right half-space to the half-space bounded by $S'$ and the $y$-axis, and add $\psi(A)$ and $\psi(B)$ to the tesselation. We iterate this procedure to cover the right half-space. Let $\mathcal V^0$ be the union of all translates of $A$ and $\mathcal V^1$ be the union of all translates of $B$. To cover the left half-space we first reflect the picture in the $y$-axis (using the isometry $\phi$). Now define
	\begin{eqnarray*}\mathcal W^0 & = & \left(\mathcal V^0\right)\cup \phi(\mathcal V^0), \\ \mathcal W^1 & = &\left(\mathcal V^1\setminus \{B_1\}\right)\cup \phi\left(\mathcal V^1\setminus \{B_1\}\right) \cup \{B_1\cup \phi(B_1)\}.
	\end{eqnarray*}
	Finally, we define $\mathcal{U}^i=\setcon{W\cap X}{W\in\mathcal{W}^i}$ for $i=0,1$.
	
	It remains to prove that $A$, $B$ and $B_1\cup \phi(B_1)$ have at most quadratic growth. Since $A$ and $B_1\cup \phi(B_1)$ is contained in the $3r$-neighbourhood of the geodesic $x=0$, they are quasi-isometric to $\Z$ and have linear growth. Now, $B$ is contained in the $r+4\delta$-neighbourhood of the union of the geodesic $x=0$ with the geodesics corresponding to upper right quadrants of circles centred at $(0,0)$ of radius $c_k=x^k$ for each $k\in\Z$. The intersection points $(0,x^k)$ are $\log(x)$ apart in the hyperbolic metric, so $B$ is quasi-isometric to the comb $C=\Z\times(\N\cup\{0\})$ with the metric $d_C((a_1,b_1),(a_2,b_2))=\abs{a_1-a_2}+b_1+b_2$. Hence, $B$ has (exactly) quadratic growth.
	It follows that $U$ is quasi-isometric to $C$ when $U\in\mathcal U^0$ and to $\Z$ when $U\in\mathcal U^1$ with quasi-isometry constants which are independent of the choice of $U$. Thus $\mathcal U\in \mathcal M_{\poly(2)}.$
\end{proof}
\begin{corollary}[{Upper bound (\ref{eq:Rank1dimformulaIntro}) in Theorem~\ref{thm:asdimseformula}}]\label{cor:asdimseupperbd}  
	Let $X = \prod_{i=1}^k X_i \times (T_3)^\ell$ where the $X_i$ are real hyperbolic spaces, $T_3$ is the infinite $3$-regular tree. Then
	\[
		\asdim_{\poly}(X)\leq \cork(X)+ \ell.
	\]
\end{corollary}
\begin{proof}
	For any $d$, there is a quasi-isometric embedding $\HH^d \to (\HH^2)^{d-1}$ (see~\cite{BradyFarb-98-filling} and the proof of Theorem \ref{prop:hypplanecover} below).
	 By Proposition~\ref{prop:products} and the trivial bound $\asdim_{\poly}\leq\asdim$ we have 
	 \begingroup
	 \allowdisplaybreaks
	\begin{align*}
		\asdim_{\poly}(X) & \leq  \sum_{i=1}^k \asdim_{\poly}(X_i) +  
		 \sum_{j=1}^\ell\asdim_{\poly}(T_3)\\
		& \leq  \sum_{i=1}^k (\dim X_i-1)\asdim_{\poly}(\HH^2) + \sum_{j=1}^\ell \asdim(T_3) \\
		& = \sum_{i=1}^k (\dim X_i-1)+ \sum_{j=1}^\ell \asdim(T_3)\\
		& = \cork(X) + \ell.\qedhere
	\end{align*}
	\endgroup
\end{proof}

\subsection{A tighter bound on the growth of covers of $\HH^k$}
Since we have $\asdim(\HH^k)=k$ but $\asdim_{\se}(\HH^k)=k-1$ it is interesting to consider $\asdim_{\poly(d)}(\HH^k)$ for different $d$ to see when the value drops down.
It is fairly straightforward to show that $\asdim_{\poly(1)}(\HH^2)=2$.  
For general $k$, using a quasi-isometric embedding of $\HH^k\to(\HH^2)^{k-1}$ and na\"ively applying the product formula (Proposition \ref{prop:products}) to Proposition~\ref{prop:tesselationH2} we can deduce that $\asdim_{\poly(2(k-1))}(\HH^k)\leq k-1$, but with more care we actually obtain the upper bound  $\asdim_{\poly(k)}(\HH^k)\leq k-1$ from Theorem \ref{prop:hypplanecover}.

We begin with a definition.

\begin{definition} Given $D>0$, a $D$-\textbf{comb} in $\HH^2$ is a union $\gamma\cup \bigcup_{i\in\Z}\gamma_i$, where $\gamma:\R\to\HH^2$ is a geodesic line and each $\gamma_i$ is the intersection of the geodesic line which intersects $\gamma$ orthogonally at the point $\gamma(iD)$ with the closure of a half-space $h$ of $\HH^2\setminus \gamma$ chosen independently of $i$.
\end{definition}
As we would like to apply our results specifically to the cover constructed in Proposition \ref{prop:tesselationH2} we additionally assume in this section that (in the upper half-space model) all of the $\gamma_i$ converge to boundary points on the $x$-axis.

\begin{lemma}\label{lem:lingrowthcones}  Whenever $C$ is a $D$-comb and $C_0$ is a connected component of $C\cap \R\times (0,\exp(a)]$, for all $c\in \R$ we have
	\[
	\left|C_0\cap (\R\times \{\exp(c)\})\right|\leq 3+\frac{2\log2}{D}+\frac{2(a_0-c)}{D}
	\]
	where $a_0=\sup\setcon{b\leq a}{C_0\cap \R\times\{\exp(b)\}\neq\emptyset}$.
\end{lemma}
\begin{proof} If $\gamma$ is a translate of the $y$-axis then this is obvious since the number of $\gamma_i\subset C_0$ which $\R\times\{\exp(c)\}$ can intersect grows linearly in $(a_0-c)$ at rate $1/D$, so assume that the boundary of $\gamma$ lies on the $x$-axis.
	
	If necessary, translate $C$ so that the point on $C$ with maximal $y$-coordinate lies on the $y$-axis. 
	
	Note that, for each $b$, $\gamma$ intersects $\R\times \{\exp(b)\}$ at most twice and each $\gamma_i$ intersects $\R\times \{\exp(b)\}$ at most once.
	Moreover, if $\gamma_i$ and $\gamma_j$ both intersect $\R\times \{\exp(b)\}$ once, then so does $\gamma_k$ for all $i\leq k \leq j$.
	
	We have $\left|C_0\cap (\R\times \exp(a_0))\right|=1$. Denote this point by $\gamma(r)$. Up to reflecting in the $y$-axis, we may assume that $\gamma(r)=(R\cos(\theta),R\sin(\theta))$ where $\theta\in[0,\pi/2)$. Moreover, $\gamma(r+s)=(R\cos(\theta+\theta_s),R\sin(\theta+\theta_s))$ where $\cos(\theta_s)^{-1}=\cosh(s)$.
	If $\gamma_i\cap \R\times \{\exp(c)\}\neq \emptyset$, then the $y$-coordinate of $\gamma(iD)$ is at least $\exp(c)$. Set $s=|iD-r|$, we have
	\begin{equation}\label{eq:addthetas}
		\frac{R\cos(\theta+\theta_s)}{R\cos(\theta)}\geq \frac{\exp(c)}{\exp(a_0)}=\exp(c-a_0)
	\end{equation}
	Rearranging, we see that
	\[
	\frac{1}{\cosh(s)}=\cos(\theta_s)\geq \cos(\theta_s)-\tan(\theta)\sin(\theta_s) \geq \frac{\exp(c)}{\exp(a_0)}=\exp(c-a_0)
	\]
	so $|iD-r|\leq a_0-c +\log(2)$. Thus, $C_0\cap (\R\times \{\exp(c)\})$ contains at most $2$ points from $\gamma$ and at most one point from each $\gamma_i$ with $|iD-r|\leq a_0-c+\log 2$. Thus
	\begin{align*}
	\left|C_0\cap (\R\times \{\exp(c)\})\right| 
		& \leq 2 + 2\frac{(a_0-c+\log 2)}{D}+1 
		\\ & \leq 3+\frac{2\log2}{D}+\frac{2(a_0-c)}{D}.
		\qedhere
	\end{align*}
\end{proof}

\begin{proof}[Proof of Theorem~\ref{prop:hypplanecover}] We will continue to work with upper half-space models. For each $k\geq 2$, the map $\alpha_k:\HH^k\to(\HH^2)^{k-1}$ given by
	\[ \alpha_k(x_1,\ldots,x_{k-1};y) = ((x_1;y),\ldots,(x_{k-1};y))
	\]
	is a quasi-isometric embedding \cite[Proposition 4.1]{BradyFarb-98-filling}.
	
	By the proofs of Propositions \ref{prop:products} and \ref{prop:tesselationH2}, for every $r>0$, there is a decomposition $(\HH^2)^{k-1}\stackrel{(r,k-1)}{\maps}\mathcal \{U_j\}$ where each $U_j$ is contained in a set of the form $C=\prod_{i=1}^{k-1} C_i$ and each $C_i\subseteq \HH^2$ is either $B_1\cup \phi(B_1)$ or the image of $A$ or $B$ under an isometry of $\HH^2$. We claim that $\alpha_k^{-1}(C)$ has polynomial growth of degree at most $k$ uniformly.
	
	Note that since $B_1\cup \phi(B_1)$ and each isometric copy of $A$ is contained in a bounded neighbourhood of an isometric copy of $B$, it suffices to assume that each $C_i=\psi_i(B)$.
	
	Let $\underline x=(x_1,\ldots,x_{k-1};y)\in \HH^k$. The ball of radius $r$ around $\underline x$ is contained in $\R^{k-1}\times[y\exp(-r),y\exp(r)]$. Each $C_i$ is a uniform neighbourhood of a $D$-comb $C'_i$ for some fixed $D$. Let $C'_i$ be a connected component of $C_i\cap \R\times[y\exp(-r),y\exp(r)]$ and define $C'=\prod_{i=1}^{k-1}C'_i$. We can bound the volume of a bounded neighbourhood of $\alpha_k^{-1}(C')\cap B(\underline{x},r)$ $\alpha_k^{-1}(C)$ (up to uniform multiplicative error) by the number of points in 
	\[
	S_r=\alpha_k^{-1}(C')\cap \bigcup_{j=-r}^r\R^{k-1}\times\{y\exp(j r)\}.
	\]
	Applying Lemma \ref{lem:lingrowthcones}, we see that
	\begin{align*}
		|S_r|  \leq  \sum_{j=-r}^r \left(3+\frac{2\log2}{D}+\frac{2(r-j)}{D}\right)^{k-1}  
		 \leq  \left(\frac2D\right)^{k-1}\sum_{l=0}^{2r} (l+K')^{k-1}
	\end{align*}
	for some constant $K'=K'(D)$. Therefore $|S_r|\leq P_k(r)$ where $P_k$ is a polynomial of degree $k$. We use the same method as in the proof of Proposition \ref{prop:tesselationH2} to pass the same result to a maximal $1$-separated subset of $\HH^k$.
\end{proof}

\section{Uniform asymptotic dimension with subexponential fibres}
\label{sec:unif-asdimse-hyp-plane}

The goal of this section is the following:

\begin{varthm}[Theorem \ref{thm:StrongPolyDimH2}] For every $d\in\N$, $\uasdim_{\poly(d)}(\HH^2)=2$.
\end{varthm}
This implies Theorem~\ref{thm:regTH2} since for any $Y$ of polynomial growth, for $d$ large enough $\uasdim_{\poly(d)}(T_3 \times Y) \leq 1$, thus there does not exist a regular map $\HH^2 \to T_3\times Y$.

Before commencing the proof, we remark that the cover constructed in the proof of Theorem \ref{prop:hypplanecover} demonstrates that $\uasdim_{\poly}(\HH^2)=1$. Therefore, Theorem \ref{thm:StrongPolyDimH2} is sharp. This is because when one takes a $D$-comb in the cover and inductively adds its neighbours finitely many times, one still has a set of polynomial growth (of higher degree).

There are three major steps to the proof of Theorem~\ref{thm:StrongPolyDimH2}, and their proofs are independent.
The main idea is that any $2$-coloured cover of $\HH^2$ must in some sense look like that constructed in Section~\ref{sec:upper-bounds-hyp-plane}.
First, we note that coarsely connected, subexponentially growing subsets of $\HH^2$ are quasi-convex quasi-trees.
Second, we show that any two-coloured partition of $\HH^2$ can be refined to give another two-coloured partition so that all sets in the partition are quasi-convex quasi-trees with at least $2$ boundary points.
Third, we find lower bounds on growth of neighbourhoods of subsets of this cover.

\subsection*{Step 1: Coarsely connected subsets of hyperbolic spaces with subexponential growth are quasi-convex quasi-trees}\ 

This essentially follows from~\cite[Proposition 2.8]{BuSh-07-hyp-tree-bound-below}; we include a proof for completeness.

\begin{lemma}\label{lem:qconv} Suppose $U\subseteq \HH^d$ is $r$-connected and has subexponential growth. Then $U$ is a quasi-convex quasi-tree. Specifically, there is some $s=s(r)$ such that for every sequence $x=x_0, x_1,\ldots, x_m=y$ in $U$ with $d(x_i,x_{i+1})\leq r$ for all $0\leq i < m$, and every point $p$ lying on the geodesic $[x,y]$ then some $x_i$ is within distance $s$ of $p$.
\end{lemma}
\begin{proof}
	We will work in the Poincar\'e ball model of $\HH^d$.
	
	Let $x=x_0,\ldots, x_m=y$ be a sequence in $U$ with $d(x_i,x_{i+1}) \leq r$ for each $i$.  (As $U$ is $r$-connected, such sequences exist for any $x,y \in U$.)  Suppose that there is a point $p\in[x,y]$ such that $d(x_i,p)> s$ for every $i$, where $s$ will be chosen later.
	
	Build a continuous path $\gamma$ from $x$ to $y$ by concatenating geodesics $[x_i,x_{i+1}]$. We automatically have $d(\gamma,p)\geq s-r/2$. For each $l\geq 1$ define $\gamma_l=\gamma \cap (B(p, (l+1)s/4)\setminus B(p, ls/4)$.
	
	Define a function $F:\HH^d\setminus\{p\}\to[0,\pi]$ where $F(z)$ is the angle between the geodesics $[x,p]$ and $[z,p]$.
	
	Let $m_l$ be the Lebesgue measure of the (possibly empty) union of closed intervals $F(\gamma_l)$.
	
	\begin{claim*} For each $k>1$ we have $\length(\gamma\cap B(p,ks/4))\geq \sum_{l=1}^{k-1} \frac{m_l}{4}e^{ls/4}$, assuming $s\geq 4$.
	\end{claim*}
	\begin{proof}[Proof of Claim] Fix a point $o\in\HH^2$ and a geodesic ray $\ell$ in $\HH^2$ starting at $o$. Define the map $\pi:\HH^d\to\HH^2$ which maps each point $z\in\HH^d$ to the unique point $y \in \HH^2$ such that $d(p,z)=d(o,y)$ and the angle between $\ell$ and $[o,y]$ measured clockwise is $F(z)$. 
		
		Next for a fixed $l$ consider the straight line retraction $\pi_l:\HH^2\to B(o, ls/4)$. The maps $\pi$ and $\pi_l$ are $1$-Lipschitz so do not increase arc-length.
		
		Let $\gamma'$ be a connected component of $\gamma_l$ with $F(\gamma')=[a,b]$. Now, $\pi_l\circ\pi(\gamma')$ is an arc of radius $ls/4$ and angle $b-a$ centred at $o$. Hence
		$\length(\gamma')
		\geq \length(\pi_l\circ\pi(\gamma'))
		=(b-a)\sinh(ls/4)\geq \frac{b-a}{4}e^{ls/4}$, 
		as $l,s/4\geq 1$.
	\end{proof}
	
	For some $l\geq 1$, $m_l\geq 2^{-l}$, since $\bigcup_{l\geq 1} F(\gamma_l)= F(\gamma)=[0,\pi]$. Hence by the claim, for this $l$,
	\[
	\length(\gamma\cap B(p,(l+1)s/4)) \geq 2^{-l-2}e^{ls/4}.
	\]
	It follows that $U\cap B(p,(l+1)s/4+r/2)$ contains at least $\frac1r 2^{-l-2}e^{ls/4}$ elements of $\{x_i\}$. Thus, assuming $s$ is large enough, we have 
	\begin{equation}\label{eq:exp}
		\abs{U\cap B(p,(l+3)s/4)} \geq \frac{1}{s} 2^{-l-2}e^{ls/4} \geq \exp(1/8)^{(l+3)s/4}.
	\end{equation}
	
	As $U$ has subexponential growth and $l\geq 1$, (\ref{eq:exp}) fails for all large enough $s$, which completes the proof.
\end{proof}
\begin{remark}
	Using the Bonk--Schramm embedding theorem, the above result can be generalised to any hyperbolic space with bounded growth at some scale (cf.\ \cite[Theorem 1.1]{BS-00-gro-hyp-embed}).
\end{remark}

\subsection*{Step 2: Refining covers.}

We begin with a key observation that allows us the flexibility to change covers while maintaining control over the sets $N^m_s(V)$ (recall Notation~\ref{notation:N}).

\begin{lemma}\label{lem:changecoverugrowth} Let $\mathcal U,\mathcal V$ be covers of a metric space $X$. If there is some constant $D$ such that, for every $V\in\mathcal V$, there is some $U_V\in\mathcal U$ such that $V\subseteq N^D_D(U_V)$, then for every $m,s$ there exist $m',s'$ such that
	\[
	N^m_s(V) \subseteq N^{m'}_{s'}(U_V).
	\]
\end{lemma}
\begin{proof} 
	We recall the definition of $N^m_s(V)$ for $V\in\cV$: we have $N^0_s(V)=V$ and for $m\geq 1$, $N^m_s(V)$ is the union of all sets in $\cV$ which intersect the closed $s$-neighbourhood of $N^{m-1}_s(V)$. 
	It follows from the definitions that for any $U,U' \in \cU$ and any $a,b,c,d$, if $N_b^a(U)$ and $N_d^c(U')$ intersect, then $N_d^c(U') \subset N_{\max\{b,d\}}^{a+2c+1}(U)$.

	Let $s'=\max\{s,D\}$.  We show that by induction on $m$, for any $m$, $N_s^m(V) \subseteq N_{s'}^{D+(2D+1)m}(U_V)$.
	The case $m=0$ is given. 
	Assume the claim is true for $m-1$.
	If $V' \in \cV$ intersects the closed $s$-neighbourhood of $N_s^{m-1}(V)$ then by induction $N_{s'}^{D+(2D+1)(m-1)}(U_V)$ and $N_D^D(U_{V'})$ intersect.
	Hence 
	\[
		V' \subset N_D^D(U_{V'}) \subset N_{s'}^{D+(2D+1)m}(U_V)
	\]
	and the claim follows. 
\end{proof}

Fix a regular tesselation $X^2$ of $\HH^2$, which we consider as a $2$-complex equipped with a metric $d_{\HH}$ which is some multiple of the metric on $\HH^2$ chosen such that $d_{\HH}(C,C') \geq 1$ for all pairs of disjoint closed cells in $X$ ($0$-, $1$- or $2$-dimensional).

Denote the $1$-skeleton of $X^2$ by $X$ and the shortest path metric on $X$ by $d_X$. Choose $M$ so that $\diam_{\HH^2}(C)\leq M$ and $\diam_X{\partial C}\leq M$ for every closed $2$-cell $C$ in $X^2$, where $\partial C \subset X$ denotes the boundary $C$. Fix some $r>10M$.

The main goal of Step 2 is to show that any cover of $\HH^2$ of multiplicity $2$ by quasi-convex quasi-trees can be `tamed' by e.g.\ discarding sets in the cover entirely enclosed by other sets.  More precisely:

\begin{proposition} Let $X\stackrel{(r,1)}{\maps} \cU$ be a decomposition of $X$ by $K$-quasi-convex $K$-quasi-trees. There exist constants $K', D$ and another decomposition $X \stackrel{(r,1)}{\maps} \cV$ of $X$ by $K'$-quasi-convex $K'$-quasi-trees with at least two boundary points, such that for every $V \in \cV$ there exists $U_V \in \cU$ with $V \subseteq N_D^D(U_V)$.
\end{proposition}

We let $\mathfrak{P}=\mathfrak{P}_0\sqcup\mathfrak{P}_1$ be a partition of $VX$ such that (with respect to the shortest path metric on $X$) elements of each $\mathfrak P_i$ are $r$-connected $K$-quasi-convex $K$-quasi-trees, and any two distinct subsets of $\mathfrak P_i$ are at distance at least $r+1$. We will refer to elements of $\mathfrak P_i$ as ``pieces'' of ``colour'' $i$.

\noindent\textbf{Refinement 1:}
\nopagebreak

Let $P\in\mathfrak P_0$ and denote by $P_r$ a subset of $X$ obtained by adding geodesics connecting every pair of vertices in $P$ at $X$-distance at most $r$. Then define $\overline P_r$ to be the sub-complex of $X^2$ obtained from $P_r$ by adding every closed $2$-cell whose vertex set is in $P$.

Define $\mathfrak{Q}_0=\setcon{\overline P_r }{P\in\mathfrak P_0}$ and set $\mathfrak{Q}_1$ to be the collection of connected components of $X^2\setminus\bigcup_{Q\in\mathfrak{Q}_0}Q$. Observe that the $1$-skeletons of pieces of $\mathfrak{Q}_0$ are connected subgraphs of the $1$-skeleton of $X$, while for every piece $Q$ of $\mathfrak{Q}_1$, the topological closure of $Q$ in $\HH^2$, $\overline{Q}$, is the union of some closed $2$-cells of $X^2$ such that at least one vertex in the boundary of each is in $\mathfrak{P}_1$. 
Since $d_X(P,P')\geq r+1$ for any distinct $P,P'\in \mathfrak{P}_0$, we have $d_X(P_r,P'_r)\geq 1$. Therefore, by assumption, $d_{\HH}(P_r,P'_r)\geq 1$.

\begin{lemma}\label{lem:reduc1}
	For every piece $Q$ of $\mathfrak{Q}_1$, there exists a unique piece $P$ of $\mathfrak{P}_1$ such that
	\[\overline{Q}\cap P\neq \emptyset.\] Moreover, $\overline{Q}$ is contained in the $2M$-neighbourhood of $P$ (for the metric on $\HH^2$), where $M$ is the diameter of a 2-cell of $X$.
\end{lemma}
\begin{proof}
	Since $\mathfrak{P}_0\sqcup\mathfrak{P}_1$ is a partition of $VX$, the existence of $P$ is obvious. Suppose towards a contradiction that $\overline{Q}$ contains two vertices $x$ and $y$ from two different pieces $P,P'$ of $\mathfrak{P}_1$. Since $\overline{Q}$ is connected, there exists a sequence of $2$-cells $C_1,\ldots,C_m$ of $\overline{Q}$ such that two consecutive cells have non-trivial intersection, $x\in C_1$, and $y\in C_m$. By assumption $d_X(P,P')\geq r+1$ which is at least 10 times the $d_X$-diameter of a $2$-cell. Let $C_j$ be the first cell which contains no vertices in $P$. Since $\diam(C_i)\leq r/10$ and $d_X(C_i,P)\leq r/10$, no vertex of $C_i$ is contained in any piece in $\mathfrak P_1$. Therefore, by construction, $C_i$ is contained in some piece in $\mathfrak{Q}_0$ which is a contradiction.
	
	The same argument shows that any closed $2$-cell in $\overline Q$ contains a vertex in $P$. Hence $\overline{Q}$ is contained in the $M$-neighbourhood of $P$.
\end{proof}

It follows from Lemma \ref{lem:reduc1} that every piece of $\mathfrak{Q}$ is contained in a (uniformly) bounded neighbourhood of a piece of $\mathfrak{P}$. Hence, increasing $K$ if necessary, pieces of $\mathfrak{Q}$ are connected $K$-quasi-convex $K$-quasi-trees.

\begin{lemma}\label{lem:domain} Let $P$ be a closed connected sub-complex of $X^2$ and let $R$ be a connected component of $\HH^2\setminus P$. If $\partial_{\infty} R\subseteq \bdry \HH^2$ has empty interior, then $|\bdry R|\leq 1$ and $R$ is contained in the union of all geodesics in $\HH^2$ with both endpoints in $P$.
\end{lemma}
\begin{proof}
	Note that $\partial R \subset \HH^2$ is contained in $P$.  Take $x \in R$.  There is at most one geodesic ray from $x$ entirely contained in $R$, for if there were two distinct such rays $\alpha,\beta$ then by the Jordan curve theorem $P$ would lie entirely on one side of $\alpha \cup \beta$, hence $\bdry R$ would contain a segment connecting $\alpha(\infty)$ and $\beta(\infty)$, contradiction.

	So if we take two different bi-infinite geodesics $\gamma_1, \gamma_2$ through $x$, then for at least one of them, say $\gamma_1$, both $\gamma_1((-\infty,0))$ and $\gamma_1((0,\infty))$ intersect $\partial R$, hence both intersect $P$.
\end{proof}

Since our pieces are $K$-quasi-convex and $\HH^2$ is hyperbolic, the union of all geodesics with endpoints in $P$ is contained in a uniformly bounded neighbourhood of $P$ by the Morse Lemma.

We introduce an order between pieces in $\mathfrak{Q}$ as follows: $Q_1<Q_2$ if $Q_1$ is contained in a connected component $R$ of the complement of $Q_2$ such that $\partial_{\infty} R$ has empty interior. We say a piece $Q_1$ is \textbf{dominated} if there exists a piece $Q_2$ such that $Q_1<Q_2$. 

\begin{lemma}\label{lem:dominatedPieces} 
	A piece is dominated if and only if it has at most one boundary point.   \end{lemma}
\begin{proof}
	If $Q_1$ is dominated by $Q_2$, with $Q_1 \subset R$ and $\bdry Q_1 \subseteq \bdry R$, then by Lemma~\ref{lem:domain} $|\bdry Q_1|\leq 1$.
	Conversely, if a piece has at most one boundary point, then by the Jordan curve theorem, there exists a unique domain in its complement whose boundary has non-empty interior, and a unique piece $Q_2$ that contains the frontier of this domain. It follows that $Q_1<Q_2$.
\end{proof}
\begin{lemma}\label{lem:risingsequence} 
	Any ordered sequence $Q_1<Q_2\ldots<Q_n$ has length $\leq L=O_{K}(1)$.  \end{lemma}
\begin{proof}
	Since $Q_i$ and $Q_{i+2}$ are disjoint, they are at distance at least $1$. Hence the closed $1$-neighbourhood of $Q_i$ must be contained in the component of the complement of $Q_{i+2}$ that contains $Q_{i+1}$. By an immediate induction, this yields the existence of a ball of radius $n/2$, which is contained in a component of the complement of $Q_n$ that contains $Q_1$. By Lemma \ref{lem:domain}, this ball is contained in the union of all geodesics with both endpoints in $P_1$. Call this union $\overline P_1$. Since  $\overline P_1$ is itself an $O(K)$-quasi-convex $O(K)$-quasi-tree, every metric ball it contains has radius at most $O(K)$. Hence $n=O(K)$, as required.
\end{proof}

\noindent\textbf{Refinement 2:} Removing dominated pieces.

Given a piece $Q$, its {\it filling} $\hat{Q}$ is the union of $Q$ with all the components of its complement whose boundary have empty interior, or equivalently, all the pieces it dominates. By Lemma \ref{lem:domain}, $\hat Q$ is contained in the union of all geodesics with endpoints in $Q$, hence it is contained within a bounded neighbourhood of $Q$ (where the bound depends only on $K$).
For every $i=0,1$, we let $\mathfrak{Q}'_i$ be the set of pieces of $\mathfrak{Q_i}$ with at least 2 boundary points. We let $\mathfrak{R}_i$ be the set of $\hat{Q}$, where $Q\in \mathfrak{Q}'_i$. This yields a new 2-coloured partition of $\HH^2$.

\begin{lemma}\label{lem:R}
	The partition $\mathfrak{R}$ is made of sets which are connected $O_{K}(1)$-quasi-convex $O_{K}(1)$-quasi-trees with at least 2 boundary points. Moreover, two distinct pieces of the same colour are at distance at least $1$.
\end{lemma}
\begin{proof} By Lemma \ref{lem:domain}, each $\hat{Q}$ is a $O_{K}(1)$-quasi-convex $O_{K}(1)$-quasi-tree. 
	Note that two pieces of same colour  are of the form $\hat{Q}_1$ and $\hat{Q}_2$, where $Q_1$ and $Q_2$ are two pieces from the same $\mathfrak{Q}'_i$.  Observe that  $Q_1$ separates its dominated pieces from $Q_2$ and its dominated pieces. Hence the distance between $\hat{Q}_1$ to $\hat{Q}_2$ is at least the distance between $Q_1$ to $Q_2$, and therefore at least $1$. 
\end{proof}

\subsection*{Step 3: Lower bounds on growth.}
We shall need the following definition.

\begin{definition}
	The {\bf comb} $C_d$ of step $d$ is a simplicial tree defined as follows. For $d=1$, it consists of a copy of the real line (seen as Cayley graph of $\Z$); $C_2$ is obtained from $C_1$ by attaching a half-line to each vertex of $C_1$: these are the {\bf hairs} of $C_2$, whose {\bf roots} are their attaching points. The embedded copy of $C_1$ will be called the {\bf base} of $C_2$.
For $d\geq 3$, we define $C_{d}$ inductively by attaching half-lines to all vertices of hairs of $C_{d-1}$ except at their roots: these are the hairs of $C_d$ (whose roots are their attaching points).
\end{definition}
By construction, combs come with natural inclusion maps $C_{d-1}\subset C_d$.	
	It is straightforward to see that a comb $C_d$ has growth $\poly(d)$.	
	
	We shall need a general classical result of hyperbolic geometry. 
	
	\begin{proposition}\label{prop:quasiconvexStatement}
	Let $Y\subset \HH^2$ be a $K$-quasi-convex subset. There exists $n_1=n_1(K)$ such that for all $L>0$, there exists $n_2=n_2(K,L)$ such that the following holds: for any two points $x,y\in Y$ such that $d(x,y)\geq n_2$, if $Z_x$ and $Z_y$ are two $K$-quasi-convex subsets such that $d(Y,Z_x\cup Z_y)\geq n_1$, $\max\{d(Z_x,x),d(Z_y,y)\}\leq L$, then $Z_x$ and $Z_y$ are contained in distance $1$ apart half-spaces whose closures are disjoint in $\overline{\HH^2}$. 
	\end{proposition} 
	
	\begin{lemma}\label{lemma:Comb} Let ${\mathfrak R}={\mathfrak R}_1\cup{\mathfrak R}_2$ be a $2$-coloured partition of $\HH^2$ where each piece is a connected $K$-quasi-convex $K$-quasi-tree with at least $2$ boundary points. 
		There exists $M=M(K), C=C(K)$ so that for every $d\geq 1$, for every piece $T_0\in\mathfrak{R}$, there is a $C$-bi-Lipschitz embedding $C_d\to N^d_M(T_0)$. In particular, $\Gr(N^d_M(T_0))\gtrsim n^d$.
	\end{lemma}

	\begin{proof}
		We proceed by induction on $d$. The case $d=1$ follows from the fact that $T_0$ is two-ended. Clearly in order to prove the lemma, it is enough to extend a bi-Lipschitz embedding of the  $(d-1)$-comb to the $d$-comb.  This amounts to showing that given a bi-Lipschitz embedded copy of $\R$ contained in  $N^{d-1}_M(T_0)$, and a half-space $A$ delimited by it, we can grow a bi-Lipschitz $2$-comb  contained in $N^{d}_{O_{K}(1)}(T_0)$ whose base is $\gamma$, and whose hairs are in $A$. Now for every $n_1>0$, since the union of pieces that intersect the $n_1$-neighbourhood of $\gamma$ in $A$ is a $O_{n_1,K}(1)$-quasi-tree, one can find a point at distance $O_{n_1,K}(1)$ from every point $x$ of $\gamma$ which lies outside these pieces. Hence there exists a piece $T_x$ that passes at  distance $\leq L=L(n_1,K)$ from $x$ but at distance at least $n_1$ from it.

		Now take a discretization of $\gamma$ whose points are sufficiently far apart--say at distance at least $n_2$. We can now apply Proposition \ref{prop:quasiconvexStatement} with $Y=\gamma$: if $n_1$ is large enough, and $n_2$ is large enough depending on $L$, then for any two of these points $x$ and $y$, their corresponding pieces $T_x$ and $T_y$ must be contained in  distance $1$ apart half-spaces $A_x$ and $A_y$ contained in $A$ delimited by quasi-geodesics $\gamma_x$ and $\gamma_y$ whose closures in $\overline{\HH^2}$ are disjoint. 
		
		In order to extend the embedding of $\gamma$ to the $2$-comb, we simply pick for every point $x$ of our discretization, an infinite $(K,K)$-quasi-geodesic $\lambda_x$ contained in $A_x$, issued from a point at distance at most $L$ from $x$. The union of $\gamma$ and the $\lambda_x$'s form a Lipschitz embedded copy of $C_2$, and since for all $x\neq y$, $\lambda_x$ and $\lambda_y$ belong to distance $1$ apart half-spaces with disjoint closures in $\overline{\HH^2}$, the embedding is quasi-isometric. 
	\end{proof}
	
	\subsection*{Conclusion of proof}
\begin{proof}[Proof of Theorem~\ref{thm:StrongPolyDimH2}]
	Suppose $\uasdim_{\poly(d)}(\HH^2) \leq 1$.
	Then for any sufficiently large $R$, there exists a decomposition $X \stackrel{(R,1)}{\maps} \cU$ of $X$ with $\cU = \{U_{ij}\}$ such that for any $s>0, m \geq 0$, $ \{N^m_s(U_{ij}))\mid U_{ij}\in \mathcal U\} \in \poly(d)$.
Following Step 2 above, $\cU$ can be modified into a two coloured partition $\mathfrak{R}$ of $\HH^2$, so that for any $\hat{Q} \in \mathfrak{R}$ there exists $P\in \cU$ with $\hat{Q}\subset N_{O(1)}(N_D^D(P))$.
	By Step 3, this contradicts our assumption for $m, s$ large enough.	
\end{proof}

	\subsection{Obstructing quasi-isometric embeddings}
	A short modification of the proof above allows us to rule out quasi-isometric embeddings $\HH^2\to T_3 \times Y$ for $Y$ with subexponential growth.

\begin{proof}[Proof of Theorem~\ref{thm:QITH2}]
	Suppose we have a quasi-isometric embedding $\HH^2\to T_3\times Y$, for $Y$ of subexponential growth, and hence a quasi-isometric embedding $f:X\to T_3\times Y$ for $X$ the approximating space for $\HH^2$ as above.
	Let $\pi:T_3\times Y \to T_3$ be the projection onto the tree, which has asymptotic dimension $1$.
	For any $R$, by pulling back a cover of $T_3$ under $f^{-1}\circ \pi^{-1}$, one can find an $R$-separated cover of $X$ by $R$-connected subsets which have uniform subexponential growth.

	The proof of Theorem~\ref{thm:StrongPolyDimH2} proceeds as above to build embedded combs $C_n$ in $\HH^2$.  These satisfy:
	\begin{itemize}
		\item 	$\cup_n C_n$ is an increasing union of $C$-bi-Lipschitz embedded combs of step $n$ (hence of growth $r^{n+1}$);
		\item  there exists a constant $D$ such that $f$ maps $C_n$ to $B(o,Dn) \times Y$.
	\end{itemize}
	Let $x$ be a based vertex in $T_0\subset \HH^2$. 
	Note that there are at least $2^n$ branches going out of the ball $B(x,Cn)$ in $C_n$ (that is: at scale $n$, $C_n$ looks like the $T_3$).
	Moreover these branches diverge linearly: for all $k\geq 2$, in the complement of $B(x, kCn)$, any two of these branches are at distance at least $(k-1)n/C$. 

	We pick $y_1,\ldots,y_{2^n}$, each in one such branch, in the ball of radius $B(x,kCn)$, such that these points are pairwise at distance $(k-1)n/C$.
	Hence, if $k$ is large enough, only depending on the quasi-isometric constants of $f$, the points $f(y_i)$ are pairwise at distance at least $3Dn$. Recall that $f(y_i)$ are contained in $B(o,Dn) \times Y$. So we deduce that their projections to $Y$ are pairwise distinct. 
	This gives $2^n$ distinct points in a ball of radius $O(n)$ in $Y$: contradiction.  
\end{proof}

\def\cprime{$'$}

\end{document}